\documentclass[letterletter, 10 pt, journal, twoside]{ieeetran}

\input epsf

\usepackage[square,numbers]{natbib}
\usepackage{xcolor,soul,framed} 
\colorlet{shadecolor}{yellow}
\usepackage{graphicx}
\usepackage{amsmath}
\newtheorem{theorem}{\textbf{Theorem}}
\newtheorem{lemma}{\textbf{Lemma}}
\newtheorem{example}{\textbf{Example}}
\newtheorem{corollary}{\textbf{Corollary}}
\newtheorem{remark}{\textbf{Remark}}
\newtheorem{definition}{\textbf{Definition}}

\usepackage{url}

\usepackage{subfigure}
\usepackage{fancyhdr}
\usepackage{amsmath}
\usepackage{multirow}
\usepackage{amssymb }
\usepackage{color}
\usepackage{graphics} 
\usepackage{graphicx}
\usepackage{amsmath}

\newenvironment{proof}{{{\bf Proof:}}}{\hfill $\square$\par}
\usepackage{multirow} 
\usepackage{xcolor}
\usepackage{multirow}
\usepackage{setspace}
\usepackage{geometry}
\usepackage{url}
\usepackage{subfigure}
\usepackage{fancyhdr}
\usepackage{amsmath}
\usepackage{multirow}
\usepackage{amssymb }
\usepackage{color}
\usepackage{graphics} 
\usepackage{graphicx}
\usepackage{algorithm} 
\usepackage{algorithmic} 
\usepackage{subeqnarray}
\usepackage{cases}
\usepackage{threeparttable}

\begin{document}
	
	\title{
		\LARGE {Zero-Norm Distance to Controllability of Linear Systems: Complexity, Bounds, and Algorithms}
	}
	\author{Yuan Zhang, Yuanqing Xia, Yufeng Zhan, and Zhongqi Sun
		
		\thanks{This work was supported in part by the
			National Natural Science Foundation of China under Grant 62003042 and Beijing Institute of Technology Research Fund Program for Young Scholars. The authors are with School of Automation, Beijing Institute of Technology, Beijing, China. Email: {\{zhangyuan14,xia\_yuanqing, yu-feng.zhan, zhongqisun\}@bit.edu.cn}. }
	}
	
	\maketitle
	
	\begin{abstract}
		Determining the distance between a controllable system to the set of uncontrollable systems, namely, the controllability radius problem, has been extensively studied in the past. However, the opposite direction, that is, determining the `distance' between an uncontrollable system to the set of controllable systems, has seldom been considered. In this paper, we address this problem by defining the notion of zero-norm distance to controllability (ZNDC) to be the smallest number of entries (parameters) in the system matrices that need to be perturbed to make the original system controllable. We show genericity exists in this problem, so that other matrix norms (such as the $2$-norm or the Frobenius norm) adopted in this notion are nonsense. For ZNDC, we show it is NP-hard to compute, even when only the state matrix can be perturbed. We then provide some nontrivial lower and upper bounds. For its computation, we provide two heuristic algorithms. The first one is by transforming the ZNDC into a problem of structural controllability of linearly parameterized systems, and then greedily selecting the candidate links according to a suitable objective function. The second one is based on the weighted $l_1$-norm relaxation and the convex-concave procedure, which is tailored for ZNDC when additional structural constraints are involved in the perturbed  parameters. Finally, we examine the performance of our proposed algorithms in several typical uncontrollable networks in multi-agent systems.
	\end{abstract}
	\begin{keywords}
		Network controllability, sparse perturbations, optimization, convex relaxation, multi-agent systems
	\end{keywords}

	\section{Introduction}
	The past decades have witnessed a renewed interest in the controllability and observability of large-scale control systems \cite{Y.Y.2011Controllability,P.Fa2014Controllability}.  Many of the real-world systems, such as social networks, transportation networks, and power networks, could be modeled or simplified as large-scale linear systems or networks of linear systems \cite{barabasi1999emergence}. Their controllability/observability, which depends on the network topology and subsystem dynamics in a complicated way \cite{L.Wa2016Controllability,Y_Zhang_2016,Composability,zhang2021structural}, is fundamental to many other system performances, including stabilization, attack detection, and secure estimation \cite{fawzi2014secure}. Considerable achievements have been made on the structure design and robustness analysis of network systems concerning controllability and observability \cite{A.Ol2014Minimal,Y_Zhang_2017,Chen2018Minimal,lou2018toward,becker2020network}.

Robustness of the controllability of a linear system against perturbations on its parameters is conventionally measured by the so-called controllability radius (CR). This notion was first proposed by Paige \cite{paige1981properties}, defined as the smallest distance (in terms of some matrix norms, such as $2$-norm and the Frobenius norm) from a given system to an uncontrollable one.  Since its  initiation, various characterizations and algorithms for CR have been proposed. To name a few, \cite{R.E1984Between} revealed an algebraic formula for CR in the complex field, that is, CR of a system $(A,B)$ is the minimum of the smallest singular value of $[\lambda I-A,B]$ with respect to $\lambda$. This relation was extended to the real field by \cite{WickDistance1991}. It was found that the real perturbations that result in uncontrollability have a rank of either one or two.  \cite{Hu2004Real} gave a unimodal formula that leads to an algorithm for computing the exact real CR using a grid search. Later, \cite{gu2006fast} gave an $O(n^4)$ bisection algorithm for computing the complex CR to any prescribed accuracy, with $n$ the number of state dimensions. A structured total least square and extended controllability matrix based algorithm was proposed in \cite{khare2012computing}.

All the above-mentioned publications deal with the unstructured CR, i.e., there is no constraint on the corresponding perturbations. Recently, a more practical scenario where the perturbations have some structural constraints has received increasing attention \cite{KarowStructured2009,bianchin2016observability,johnson2018structured,zhang2022real}. For example, \cite{bianchin2016observability} provided Lagrange multiplier-based characterizations and a heuristic
algorithm for the observability radius (the dual to CR), where the perturbations have a prescribed sparsity pattern. In \cite{johnson2018structured} and \cite{zhang2022real},  the structured CR is considered where the perturbations are affinely parameterized, with iterative algorithms proposed to find the local optima.  It was also proven in \cite{zhang2022real} that computing such a CR is NP-hard.

Altogether, the CR problem as well as its various variants has been extensively explored. However, a natural problem in the opposite direction of CR, i.e., determining the smallest distance from an uncontrollable system to a controllable one, has seldom been addressed. This may be because that controllability is a generic property in the sense that any uncontrollable system is `arbitrarily' close to a controllable one \cite{C.T.1974Structural} (in terms of the $2$-norm or Frobenius norm; see the argument in Section \ref{problem-formulation}). Nevertheless, if we consider the `zero-norm' measure, i.e., the number of entries that need to be perturbed to make a system controllable (note that `zero-norm' is not truly a norm), things become different.  To be specific, even when the perturbation transforming an uncontrollable system to a controllable one can be arbitrarily small (measured by a formal matrix norm), its `zero-norm' should be no less than a positive number. On the other hand, because of the genericity of controllability, if we find the zero-nonzero pattern of a perturbation that makes the original system controllable, then assigning random values to the nonzero entries of this pattern (i.e., almost every realization of this pattern) will lead to controllability with probability one. Hence, introducing the zero-norm distance to controllability (ZNDC), defined as the smallest number of entries in the system matrices that need to be perturbed to make a system controllable, is reasonable and desirable. This is particularly important in designing controllable networks, as ZNDC tells us how to transform an uncontrollable network into a controllable one by perturbing the smallest number of edge weights.




	
	In the context of structured systems, the distance to controllability is measured by the number of free entries that need to be added for achieving structural controllability \cite{Y_Zhang_2017}. It is shown this index can be computed in polynomial time \cite{Y_Zhang_2017,Chen2018Minimal}. But, additional constraints on the addable entries will make its computation NP-hard \cite{zhang2019minimal}. Note in the structured system theory, the free entries take values outside certain hypersurfaces, while in ZNDC, the system matrices are all numerically fixed. Hence, the notion ZNDC differs from \cite{Y_Zhang_2017}.   On the other hand, ZNDC can be regarded as a generalization of the minimal controllability problem (MCP) considered in \cite{A.Ol2014Minimal}, which seeks to find the sparsest input matrices making a system controllable. Say, if we are given a state matrix and a zero input matrix, when restricting that only entries in the input matrix can be perturbed, then the ZNDC reduces to the MCP. More recently, \cite{becker2020network} has considered perturbing a subset of edge weights to make a network system possess certain prescribed controllability metrics.

In this paper,  we characterize ZNDC in terms of its computational complexity and lower/upper bounds. Two heuristic algorithms for its computation are also given. The main contributions of this paper are three-fold:

\begin{itemize}
	\item Complexity: We show computing the ZNDC is NP-hard, even when only the state matrix can be perturbed. 
	
	\item Bounds: We give several non-trivial upper and lower bounds for ZNDC.
	
	\item Algorithms: We provide two heuristic algorithms for ZNDC. The first one is a greedy algorithm built on the structural controllability of a linearly parameterized plant, and the second one is based on the weighted $l_1$-norm relaxation and the convex-concave procedure.
	
\end{itemize}

While results on the complexity and bounds are devoted to the unstructured ZNDC, the proposed algorithms are valid when additional structural constraints are involved in the corresponding perturbations.
As complementary to CR, our results may deepen our understanding between controllability and uncontrollability.  In obtaining the above-mentioned results, we leverage the structural controllability theory of linear-parameterized plants in \cite{Morse_1976} and transform the ZNDC into a problem of structural controllability of a specifically linear-parameterized system.  

The rest is organized as follows. Section \ref{problem-formulation} gives the problem formulation. Section \ref{linear-para-sec} provides preliminaries on structural controllability of linear-parameterized plants, based on which a simplified criterion tailored for ZNDC is obtained. The next two sections characterize the computational complexity and lower/upper bounds for ZNDC. Section \ref{alg-sec} provides two
heuristic algorithms for computing ZNDC, followed by some typical examples and simulations to validate the effectiveness of the proposed algorithms in Section \ref{simu-sec}. The last section concludes this paper.

Notations: For $n\in {\mathbb N}$, let $[n]\doteq
\{1,2,...,n\}$. $I_n$ denotes the $n\times n$ identify matrix. $|\cdot|$ takes the absolute value of a scalar. For a matrix $M$, $\sigma(M)$ is the set of eigenvalues of $M$. For $M\in {\mathbb R}^{n_1\times n_2}$, $S_1\subseteq [n_1]$, and  $S_2\subseteq [n_2]$, $M(S_1,:)$ ($M(:,S_2)$) denotes the submatrix of $M$ with rows (columns) indexed by $S_1$ ($S_2$), and $M_{S_1,S_2}$ the submatrix with rows indexed by $S_1$ and columns by $S_2$. $0_{m\times n}$ ($1_{m\times n}$) denotes the $m\times n$ matrix with all entries $0$ ($1$). $M \succeq 0$ means $M$ is semi-positive definite.

	
	\section{Problem Formulation} \label{problem-formulation}
	
	Consider the following linear time-invariant system:
\begin{equation} \label{StateSpace} \dot x(t) = Ax(t) + Bu(t)\end{equation}
where $x(t)\in {{\mathbb R}^{n}},u(t)\in {{\mathbb R}^{m}}$ are respectively state vectors and input vectors,  and $A \in {{\mathbb R}^{n \times n}},B \in {{\mathbb R}^{n \times m}}$ are state and input matrices. Controllability of system (1) is the ability to steer $x(t)$ arbitrarily by choosing feasible input $u(t)$.

We consider the following problem:

\begin{equation}\tag{${\cal P}_1$} \label{prob1}  \begin{array}{l}
\mathop {\min }\limits_{\Delta A \in {{\mathbb R}^{n \times n}},\Delta B \in {{\mathbb R}^{n \times m}}} {\kern 1pt} {\kern 1pt} {\kern 1pt} {\kern 1pt} {\kern 1pt} {\left\| {\left[ {\Delta A,\Delta B} \right]} \right\|_0}\\
{\rm s.t.}{\kern 1pt} {\kern 1pt} {\kern 1pt} {\kern 1pt} {\kern 1pt} {\kern 1pt} {\kern 1pt} {\kern 1pt} {\kern 1pt} {\kern 1pt} {\kern 1pt} {\kern 1pt} {\kern 1pt} (A + \Delta A,B + \Delta B){\kern 1pt} {\kern 1pt} {\kern 1pt} {\rm{is}}{\kern 1pt} {\kern 1pt} {\kern 1pt} {\rm{controllable}}
\end{array}\end{equation}
where $||\cdot||_0$ takes the number of nonzero entries in a matrix.
To avoid the trivial case, assume that $(A,B)$ is uncontrollable. Then, \ref{prob1} aims to find the minimal number of entries of $(A,B)$ that need to be perturbed such that the resulted system is controllable. As mentioned earlier, we call the optimal value of \ref{prob1} the ZNDC (i.e., zero-norm distance to controllability) of $(A,B)$, denoted by $r_c(A,B)$. If $(A,B)$ is clear from the context, we will drop $(A,B)$ from $r_c(A,B)$.  

It is obvious that the above definition is well-defined, as $r_c\in {\mathbb N}$ and for an uncontrollable pair $(A,B)$, $r_c>0$. Besides, whenever $m>0$, there always exists feasible solutions to \ref{prob1}. For example, one can perturb all entries of $(A,B)$ (i.e., the number of perturbed entries is $n^2+nm$) and get a controllable pair $(A+\Delta A, B+\Delta B)$.

The definition of ZNCD is motivated by the {\emph{distance to uncontrollability}}, i.e., the CR of $(A,B)$, which is defined as
{\small	$$ r_{\bar c}^{|| \bullet ||} = \min \{ \left\| {\left[ {\Delta A,\Delta B} \right]} \right\|: (A + \Delta A,B + \Delta B){\kern 1pt} {\kern 1pt} {\kern 1pt} {\rm{is}}{\kern 1pt} {\kern 1pt} {\kern 1pt} {\rm{uncontrollable}}\} $$}where the matrix norm $||\bullet ||$ is the 2-norm or Frobenius norm, and $\Delta A, \Delta B$ can be in the real or complex field. In contrast to CR, if the zero-norm in \ref{prob1} is replaced with the 2-norm or Frobenius norm, then the corresponding definition becomes meaningless. That is because the distance from uncontrollability to controllability is almost zero with respect to the 2-norm or the Frobenius norm. To see this, suppose there is a $[\Delta A,\Delta B]$ such that $(A+\Delta A,B+\Delta B)$ is controllable. Consider a new pair constructed as
$$\lambda [A+\Delta A,B+\Delta B]\text{+(1-}\lambda \text{)}\left[ A,B \right]=\left[ A+\lambda \Delta A,B+\lambda \Delta B \right]$$
where $\lambda\in {\mathbb R}$. It turns out that when $\lambda=1$, the constructed pair $(A+\lambda \Delta A, B+\lambda \Delta B)$ is controllable. Hence, there is an $n\times n$ submatrix of the controllability matrix of $(A+\lambda \Delta A,B+\lambda \Delta B)$, given by ${\cal C}(A+\lambda \Delta A,B+\lambda \Delta B)$, expressed as {\footnotesize $$\left[ B+\lambda \Delta B, (A+\lambda \Delta A)(B+\lambda \Delta B),\cdots ,{{(A+\lambda \Delta A)}^{n-1}}(B+\lambda \Delta B) \right],$$}whose determinant is a nonzero polynomial of $\lambda$. Consequently, for arbitrarily small $\varepsilon >0$, there exists $0<\lambda <\varepsilon $ such that the aforementioned polynomial is not zero, i.e., the pair  $(A+\lambda \Delta A,B+\lambda \Delta B)$ is controllable. Notice that the 2-norm or Frobenius-norm of $[\lambda \Delta A,\lambda \Delta B]$ can arbitrarily approach zero as $\lambda$ approaches zero.

Moreover, in physical plants/networks, the entries of the state matrix $A$ or the input matrix $B$ that can be perturbed may have structural constraints. Particularly, consider that the system matrices are parameterized by the vector $\theta=[\theta_1,...,\theta_l]^{\intercal}$ as
\begin{equation} \label{affine}
A(\theta)=A+\sum \nolimits_{i=1}^l\theta_i A_i, B(\theta)=B+\sum \nolimits_{i=1}^l\theta_i B_i.
\end{equation}
The affine parameterization (\ref{affine}) is common in the literature for describing how system matrices are affected by the parameters \cite{Morse_1976,Anderson_1982,zhang2019structural,KarowStructured2009}. Without losing generality, consider the following optimization problem:

\begin{equation}\tag{${\cal P}_2$} \label{prob2}  \begin{array}{l}
\mathop {\min }\limits_{\theta\in {\mathbb R}^l} {\kern 1pt} {\kern 1pt} {\kern 1pt} {\kern 1pt} {\kern 1pt} {\left\|\theta\right\|_0}\\
{\rm s.t.}{\kern 1pt} {\kern 1pt} {\kern 1pt} {\kern 1pt} {\kern 1pt} {\kern 1pt} {\kern 1pt} {\kern 1pt} {\kern 1pt} {\kern 1pt} {\kern 1pt} {\kern 1pt} {\kern 1pt} (A(\theta), B(\theta)) \ {\rm in} \ (\ref{affine}) \ {\rm{is}}{\kern 1pt} {\kern 1pt} {\kern 1pt} {\rm{controllable}}
\end{array}\end{equation}
In other words, \ref{prob2} intends to find the smallest number of parameters in $\theta$ that need to be perturbed such that $(A,B)$ becomes controllable. Here, when $\theta=0_{l\times 1}$, we get the nominal $(A,B)$ of the considered system.


In the sequel, we will analyze the computation complexity, lower/upper bounds, and provide some heuristic algorithms for the above two problems.



\section{Structural controllability of a linear-parameterized plant} \label{linear-para-sec}
In this section, we recall the structural controllability of a linear-parameterized plant and reformulate the existing criteria in \cite{Morse_1976} and \cite{zhang2019structural} to fit the analysis of the addressed problem in this paper.

In \cite{Morse_1976}, controllability of a linear-parameterized pair $(A,B)$ is concerned, which is modeled as
\begin{equation} \label{Linear Parameterization} A=A_0+\sum\nolimits_{i=1}^k g_is_ih_{1i}^\intercal, B=B_0+\sum\nolimits_{i=1}^k g_is_ih^\intercal_{2i},\end{equation}
where $g_i,h_{1i}\in {\mathbb R}^n$, $h_{2i}\in {\mathbb R}^m$, $\{s_1,\cdots,s_k\}$ are real free parameters. Denote $g=[g_1,...,g_k]$, $h_1=[h_{11},...,h_{1k}]^{\intercal}$, and $h_2=[h_{21},...,h_{2k}]^{\intercal}$.

\begin{definition} \label{definition of structural controllability}
	System (1) with $(A,B)$ parameterized as (\ref{Linear Parameterization}) is said to be structurally controllable, if there exists a set of real values $\{s_1, \cdots ,s_k\}$ such that the corresponding numerically specified system (1) is controllable.
\end{definition}


To present the result of \cite{Morse_1976}, introduce two transfer functions as follows
$$ \begin{array}{l} G_{1}(\lambda)=[h_{11},...,h_{1k}]^{\intercal}(\lambda I- A_0)^{-1}[g_1,...,g_k],\\ G_{2}(\lambda)=[h_{11},...,h_{1k}]^{\intercal}(\lambda I- A_0)^{-1}B_0+[h_{21},...,h_{2k}]^{\intercal}. \end{array}$$
Associated with $G_{1}(\lambda)$ and $G_{2}(\lambda)$, an auxiliary digraph is constructed as ${\cal G}_{\rm d}=(V_z\cup V_u, E_{zz}\cup E_{uz})$, where $V_z=\{z_1,...,z_k\}$, $V_u=\{u_1,...,u_m\}$, $E_{zz}=\{(z_i,z_j): [G_{1}(\lambda)]_{ji}\ne 0\}$ and $E_{uz}=\{(u_i,z_j): [G_{2}(\lambda)]_{ji}\ne 0\}$. A vertex $z_i\in V_z$ is input-reachable if there is a path from a vertex $u_i\in V_u$ to $z_i$ in ${\cal G}_{\rm d}$. A cycle is said to be input-reachable, if at least one of its vertices is input-reachable.
The following lemma gives a necessary and sufficient condition for $(A,B)$ in (\ref{Linear Parameterization}) to be structurally controllable.

\begin{lemma}\cite{Morse_1976,zhang2019structural} \label{Theorem of Linear Parameterization}
	$(A,B)$ in (\ref{Linear Parameterization}) is structurally controllable, if and only if
	
	(a) Every cycle of ${\cal G}_{\rm d}$ is input-reachable;
	
	(b) For each $S\subseteq \{1,...,k\}$, ${\rm rank} \left[
	\begin{array}{ccc}
	\lambda I- A_0 & -B_0 & g(:,S) \\
	h_{1}([k]\backslash S,:) & h_2([k]\backslash S,:) & 0 \\
	\end{array}
	\right]\ge n
	$, $\forall \lambda \in \sigma(A_0)$.
\end{lemma}

To leverage Lemma {\ref{Theorem of Linear Parameterization}} for \ref{prob1}, we first introduce some notions. A pattern matrix ${\cal M}$ is a matrix with entries from $0$ and $*$. We use $\{0,*\}^{n_1\times n_2}$ to denote the set of all pattern matrices with dimension $n_1\times n_2$. For ${\cal M}\in \{0,*\}^{n_1\times n_2}$, let $[{\cal M}]=\{M\in {\mathbb R}^{n_1\times n_2}: M_{ij}=0 \ {\rm if}\ {\cal M}_{ij}=0\}$. Let $\bar I_{k}$ denote the $k$ dimensional diagonal matrix with diagonal entries being $*$.  The generic rank of ${\cal M}$ (given by ${\rm grank} {\cal M}$) is the maximum rank an element in $[{\cal M}]$ can achieve. A matroid is a pair ${\cal M}=(V,{\cal I})$ of a set $V$ and a collection ${\cal I}$ of subsets of $V$ satisfying: (1) $\emptyset \in {\cal I}$; (2) $I \subseteq J \in {\cal I}$ implies $I\in {\cal I}$; (3) for any $I,J\in {\cal I}$, $|I|<|J|$, there is some $v\in J\backslash I$ such that  $I\cup \{v\}\in {\cal I}$.  Here, an element of ${\cal I}$ is called an independent set.  For two matroids ${\cal M}_1=(V,{\cal I}_1)$ and ${\cal M}_2=(V, {\cal I}_2)$, the cardinality of their intersection is the maximum size of a common independent set, i.e., $\max \{|I|: I\in {\cal I}_1, I\in {\cal I}_2\}$, which value is denoted by $\rho({\cal M}_1\cap {\cal M}_2)$ and can be computed in polynomial time \cite{Murota_Book}.  Given a $p\times q$ matrix $W$, the matroid formed by columns of $W$ is defined as ${\cal M}(W)=([q], {\cal I}_W)$, with ${\cal I}_W=\{I\subseteq [q]: {\rm rank} W(:,I)=|I|\}$ (if $M$ contains a pattern submatrix, `${\rm rank}$' shall be replaced with `${\rm grank}$'.)

Rewrite $(A,B)$ with the perturbation $[\Delta A, \Delta B]$ as
\begin{equation}\label{state space perturbation} [A + \Delta A,B + \Delta B] = [A,B] + {I_n}[\Delta A,\Delta B]\left[ {\begin{array}{*{20}{c}}
	{{I_n}}&{0_{n\times m}}\\
	{0_{m\times n}}&{{I_m}}
	\end{array}} .\right]\end{equation}Let ${\cal A}\in \{0,*\}^{n\times n}$ (${\cal B}\in \{0,*\}^{n\times m}$) be the pattern matrix specifying the sparsity pattern of $\Delta A$ ($\Delta B$), that is, $\Delta A\in [{\cal A}]$ ($\Delta B\in [{\cal B}]$). Based on ({\ref{state space perturbation}}), introduce two transfer functions $${G_{zx}}(\lambda ) = \left[ \begin{array}{l}



{I_n}\\
{0_{m \times n}}
\end{array} \right]{(\lambda I - A)^{ - 1}}{I_n} = \left[ {\begin{array}{*{20}{c}}
	{{{(\lambda I - A)}^{ - 1}}}\\
	{{0_{m \times n}}}
	\end{array}} \right],$$ $${G_{zu}}(\lambda ) = \left[ \begin{array}{l}
{I_n}\\
{0}
\end{array} \right]{(\lambda I - A)^{ - 1}}B + \left[ \begin{array}{l}
{0}\\
{I_m}
\end{array} \right]{\rm{ = }}\left[ \begin{array}{l}
{(\lambda I - A)^{ - 1}}B\\
\ \ \ {I_m}
\end{array} \right].$$

\begin{figure}
		\centering
		\includegraphics[width=2.8in]{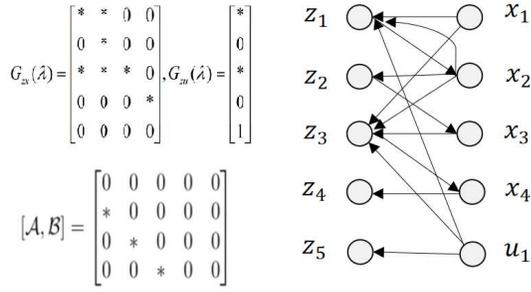}\\
		\caption{Example of the ACG. For simplicity, we only give the sparsity patterns of $G_{zx}(\lambda)$ and $G_{zu}(\lambda)$.}\label{auc1}
	\end{figure}

Construct the auxiliary connection digraph (ACG) ${\cal G}_{\rm auc}=(V_{\rm auc},E_{\rm auc})$ associated with system ({\ref{state space perturbation}}) as follows: the vertex set is $V_{\rm auc}=V_z\cup V_x\cup V_u$, with $V_z=\{{{z}_{1}},...,{{z}_{n+m}}\}$, $V_x=\{ {{x}_{1}},...{{x}_{n}}\}$, $V_u=\{u_1,...,u_n\}$, the edge set $E_{\rm auc}=E_{xz} \cup E_{uz} \cup E_{zx}$ with $E_{xz}=\{(x_j,z_i): [{G_{zx}}(\lambda )]_{ij} \ne 0\}$, $E_{uz}=\{(u_j,z_i): [G_{zu}(\lambda)]_{ij}\ne 0\}$, $E_{zx}=\{(z_i,x_j): [{\cal A}, {\cal B}]_{ji}\ne 0\}$. See Fig. \ref{auc1} for illustration ($V_z$ could be understood as a duplication of $V_x\cup V_u$).
Similar to ${\cal G}_d$, a vertex is input-reachable in ${\cal G}_{\rm auc}$, if there is a path from $V_u$ ending at it.
We say an edge $e\in E_{xz}$ is input-reachable if either the beginning  or the ending vertex of $e$ is input-reachable. 

%
%

If there exist $\Delta A \in \left[ {{\cal A}} \right]$ and $\Delta B \in \left[ {{\cal B}} \right]$ such that $(A+\Delta A, B+\Delta B)$ is controllable, we say system (\ref{state space perturbation}) is {\emph{structurally controllable}}.
With the construction of the ACG and Lemma \ref{Theorem of Linear Parameterization}, we have the following criterion for structural controllability of system (\ref{state space perturbation}).

\begin{theorem} \label{theorem 1}  For a given pattern pair ${\cal A}\in {{\{0,*\}}^{n\times n}},{\cal B}\in {{\{0,*\}}^{n\times m}}$, the following statements are equivalent:
	
	1) System (\ref{state space perturbation}) is structurally controllable.
	
	2) (a) For each $\lambda_i \in \sigma(A)$, $$\mathop {\rm{max}} \limits_{\Delta A \in [{\cal A}],\Delta B \in [{\cal B}]} {{\rm{rank}}} {\rm{ }}[{\lambda _i}I - A - \Delta A,B + \Delta B] = n;$$
	
	(b) Each edge $e\in E_{xz}$ is input-reachable in ${\cal G}_{\rm auc}$.
	
	3) (a) For each $\lambda_i \in \sigma(A)$, the intersection of two matroids formed by columns of $\left[ \begin{matrix}
	I_n & {{\lambda }_{i}}I-A & B  \\
	\end{matrix} \right]$ and $\left[\begin{matrix} [{\cal A}, {\cal B}]^\intercal & {{\bar I_{n + m}}} \\ \end{matrix} \right]$ has cardinality $n$;
	
	(b) Every vertex $z_i\in V_z$ is input-reachable in ${\cal G}_{\rm auc}$.
\end{theorem}
	
	{\bf{Proof.}}  By rewriting $(A+\Delta A, B+\Delta B)$ as ({\ref{state space perturbation}}), Condition (a) of Statement 2) follows from \citep[Proposition 3]{zhang2019structural}, and condition (b) from \citep[Corollary 1]{zhang2019structural}. Condition (a) of Statement 3) results from \citep[Proposition 5]{zhang2019structural}. For Condition (b) of Statement 3), notice that the transfer function ${{(\lambda I-A)}^{-1}}=\frac{1}{\det (\lambda I-A)}\text{adj}(\lambda I-A)$, where ${\rm adj}(\cdot)$ denotes the adjacency matrix. It yields all diagonal entries of ${{(\lambda I-A)}^{-1}}$ are nonzeros.
This implies, the edge $(x_i,z_i)\in E_{xz}$,  for $i=1,...,n$. Hence, the input-reachability of all edges in $\{(x_1,z_1),...,(x_n,z_n)\}$ indicates that each $z_i\in V_z$ is input-reachable. On the other hand, the input-reachability of every $z_i\in V_z$ certainly leads to that all $e\in E_{xz}$ are input-reachable. Therefore, Conditions (b) in Statements 3) and 2) are equivalent.  $\hfill\blacksquare$
	


	\begin{remark}
		The difference between Theorem \ref{theorem 1} and Lemma \ref{Theorem of Linear Parameterization} results from that  $[\Delta A, \Delta B]$ in (\ref{state space perturbation}) is not diagonal as in (\ref{Linear Parameterization}).
	\end{remark}
	
	Now we give an equivalent condition of Condition (a) in Statement 2) of Theorem \ref{theorem 1}.
	
	\begin{lemma} \label{lemma-int} Condition (a)  of Statement 2) in Theorem \ref{theorem 1} is equivalent to that, for each $\lambda_i \in \sigma(A)$, there exist $J_{i}^{R}\subseteq [n]$, $J_{i}^{C}\subseteq [n+m]$ such that: (i) $|J_{i}^{C}|=|J_{i}^{R}|$;
		(ii) $\bar{J}_{i}^{R}\subseteq [n]\backslash J_{i}^{R}$, $\bar{J}_{i}^{C}\subseteq [n\text{+}m]\backslash J_{i}^{C}$, and  $|\bar{j}_{i}^{C}|=|\bar{j}_{i}^{R}|\text{=}n-|j_{i}^{C}|$, such that  $[\lambda_iI -A, B]_{J_i^R,J_i^C}$ is of full rank and ${{[{\cal A},{\cal B}]}_{\bar{j}_{i}^{R},\bar{j}_{i}^{C}}}$ is of full generic rank.
	\end{lemma}
	
	{\bf{Proof.}} This lemma is a direct result of \citep[Theorem. 4.2.8]{Murota_Book}: For a matrix $M=Q+T$, where $Q$ is a constant matrix, $T$ is a structured matrix whose nonzero entries are free parameters. Then,
	$\max\nolimits_{T'\in [T]} \text{rank}\ Q+ T'=\max \{\text{rank}{{Q}_{I,J}}+\text{grank}{{T}_{R\backslash I,C\backslash J}}|I\subseteq R,J\subseteq C\}$, where $R,C$ are respectively the sets of row indices and column indices of $M$. $\hfill\blacksquare$
	
	
	For notation simplicity, if system (\ref{state space perturbation}) is structurally controllable with the pattern pair $({\cal A},{\cal B})$, we say $({\cal A},{\cal B})$ is a feasible solution to \ref{prob1}.
	
	\section{Computational complexity}

In this section, we prove that \ref{prob1} is NP-hard, even when the perturbed entries are restricted in $A$. Note \ref{prob1} differs from the MCP in \cite{A.Ol2014Minimal} given as follows, in the sense that both the input matrix $B$ and state matrix $A$ can be perturbed.

{\bf{MCP:}} Given $A\in {\mathbb R}^{n\times n}$, find
$$\begin{array}{l}
\mathop {\min }\limits_{{\kern 1pt} {\kern 1pt} {\kern 1pt} {\kern 1pt} {\kern 1pt} B \in {{\mathbb R}^{n \times m}}} {\kern 1pt} {\kern 1pt} {\kern 1pt} {\kern 1pt} {\kern 1pt} {\left\| B \right\|_0}\\
{\rm s.t.}{\kern 1pt} {\kern 1pt} {\kern 1pt} {\kern 1pt} {\kern 1pt} {\kern 1pt} {\kern 1pt} {\kern 1pt} {\kern 1pt} {\kern 1pt} {\kern 1pt} {\kern 1pt} {\kern 1pt} (A,B){\kern 1pt} {\kern 1pt} {\kern 1pt} {\rm{is}}{\kern 1pt} {\kern 1pt} {\kern 1pt} {\rm{controllable.}}
\end{array}$$

\begin{theorem} \label{complexity-1}
	\ref{prob1} is NP-hard.
\end{theorem}

\begin{proof} We give a reduction from the MCP to \ref{prob1}. For an $A\in {\mathbb R}^{n\times n}$,  assume that $A$ has no repeated eigenvalues.
	From \cite{A.Ol2014Minimal}, the MCP of finding the sparsest $B\in {\mathbb R}^{n\times 1}$ such that $(A,B)$ is controllable is NP-hard. Let $B_0=0_{n\times 1}$, and let $x_1,...,x_n$ be $n$ linearly independent left eigenvectors of $A$,  $x_i\in {\mathbb C}^{n}$, $i=1,...,n$. Then, according to the PBH test, the MCP associated with $A$ is equivalent to finding the sparsest $B^*\in {\mathbb R}^{n\times 1}$ such that
	\begin{equation} \label{numerical-cover-condition} x_i^{\intercal} B^*\ne 0 , \forall i\in \{1,...,n\}.  \end{equation}
	Let ${\cal B}^*\in \{0,*\}^{n\times 1}$ represent the sparsity pattern of $B^*$.  We will show, $[0_{n\times n}, {\cal B}^*]$ is an optimal solution to \ref{prob1} associated with $(A, B_0)$.

	The feasibility of $[0_{n\times n}, {\cal B}_*]$ for \ref{prob1} on $(A,B_0)$ is obvious. It is shown that the optimal solution to \ref{prob1} has a sparsity at least $||{\cal B}^*||_0$. To this end, let ${\cal S}(x_i)=\{j: [x_i]_j\ne 0\}$ be the support of $x_i$, and  define ${\cal S}({\cal B}^*)$ similarly. Then, for (\ref{numerical-cover-condition}) to hold, it is necessary and sufficient that
	\begin{equation} \label{cover-condition} {\cal S}(x_i)\cap {\cal S}({\cal B}^*)\ne \emptyset, \forall i\in \{1,...,n\}.    \end{equation}
	On the other hand, notice that for any feasible solution $[{\cal A}', {\cal B}' ]$ to \ref{prob1} on $(A,B_0)$,  it must hold
	\begin{equation} \label{nece-prob1} x_i^{\intercal} [{\cal A}', {\cal B}']\ne 0, \forall i\in \{1,...,n\}.  \end{equation}
	Indeed, if the condition above is not satisfied for some $i\in \{1,...,n\}$, we have $x_i^\intercal [\lambda_i I- A - \Delta A, B_0+\Delta B]=0$, for any $\Delta A \in [{\cal A}']$, $\Delta B\in [{\cal B}']$, leading to the uncontrollability of $(A+\Delta A, B_0+\Delta B)$ (recalling $B_0=0_{n\times 1}$), where $\lambda_i$ is the eigenvalue associated with the eigenvector $x_i$. Let ${\cal S}([{\cal A}',{\cal B}'])=\{j: [{\cal A}',{\cal B}'](\{j\},:)\ne 0\}$. Again, it is seen easily that for (\ref{nece-prob1}) to hold, it is necessary that
	$${\cal S}(x_i)\cap {\cal S}([{{\cal A}'}, {{\cal B}'}])\ne \emptyset, \forall i\in \{1,...,n\}.$$
	Since ${\cal B}^*$ is the sparsest ${\cal B}^*\in \{0,*\}^{n\times 1}$ satisfying (\ref{cover-condition}), we have $||[{\cal A}', {\cal B}'] ||_0 \ge ||{\cal B}^*||_0$. Therefore, $[0_{n\times n}, {\cal B}^*]$ is an optimal solution to \ref{prob1} associated with $(A, B_0)$. Since finding ${\cal B}^*$ is NP-hard (from the NP-hardness of the MCP), and the reduction above is in polynomial time, we attain that \ref{prob1} on $(A,B_0)$ is NP-hard.
\end{proof}

Based on Theorem \ref{complexity-1}, we give the complexity of \ref{prob1} when only the state matrices can be perturbed. Such a scenario may be common for large-scale distributed systems, where the input structure is often fixed (being the dedicated input structure, i.e., each input actuates only one state variable) \cite{P.Fa2014Controllability}.
\begin{corollary}
	\ref{prob1} is NP-hard, even when the perturbed entries are restricted in $A$.
\end{corollary}

\begin{proof}
	Consider the system $(A,B)$ given in (\ref{StateSpace}). Construct a new system $(A',B')$ as
	$$A'=\left[\begin{matrix} A & B \\
	0_{m\times n} & 0_{m\times m}
	\end{matrix} \right], B'=\left[\begin{matrix}  0_{n\times m} \\
	I_{m}
	\end{matrix}	\right].$$
	Suppose $(\Delta A, \Delta B)$ is an optimal solution to \ref{prob1} with $(A,B)$. Notice that
	$$\begin{array}{c}{\rm rank}\left[\begin{matrix}
	A+\Delta A-\lambda I_n  & B+\Delta B & 0 \\
	0 & -\lambda I_m & I_m
	\end{matrix}\right] = \\
	m+ {\rm rank}\left[A+\Delta A-\lambda I, B+\Delta B \right], \forall \lambda \in {\mathbb C}. 	\end{array} $$
	Hence, ${\tiny \left[\begin{matrix}
		\Delta A & \Delta B \\
		0_{m\times n} & 0_{m\times m}
		\end{matrix} \right]}$ is also feasible for \ref{prob1} with $(A', B')$ when only $A'$ can be perturbed (denote such a problem by ${\cal P}_1'$).
	On the other hand, any feasible solution to ${\cal P}_1'$ must contain a submatrix $[\Delta A', \Delta B']$ such that
	${\rm rank}\left[A+\Delta A'-\lambda I, B+\Delta B' \right]=n$, $\forall \lambda \in {\mathbb C}$. Since $(\Delta A, \Delta B)$ is the optimal one, we attain that ${\tiny \left[\begin{matrix}
		\Delta A & \Delta B \\
		0_{m\times n} & 0_{m\times m}
		\end{matrix} \right]}$ is an optimal solution to ${\cal P}_1'$. Since \ref{prob1} is NP-hard,  it follows immediately that  ${\cal P}_1'$ is also NP-hard.
\end{proof}

\section{Upper and Lower Bounds}

Since computing ZNDC is NP-hard, before presenting the heuristic algorithms, we give two nontrivial bounds for it.

\begin{theorem} \label{bounds} (upper/lower bounds)  Given $(A,B)$, the ZNDC $r_c$ satisfies
	
	(1) $r_c\le n-{\rm rank} B$;
	
	(2) $r_c$ is no less than the optimal value of the following problem
	\begin{equation}\label{lower bound}\begin{array}{l}
	\mathop {\min }\limits_{J \subseteq \{ 1,...,n\} } {\mkern 1mu} \left| J \right|\\
	{\rm s.t.}{\kern 1pt} {\kern 1pt} {\kern 1pt} {\kern 1pt} {\kern 1pt} {\kern 1pt} (A,[B,{I(:,J)}]){\kern 1pt} {\kern 1pt} {\kern 1pt} {\rm{is}}{\kern 1pt} {\kern 1pt} {\kern 1pt} {\rm{controllable}}
	\end{array}\end{equation}
\end{theorem}

{{\bf Proof.}} We first prove the upper bound. We divide the proof into two cases.

{\bf{Case I:}} Assume ${\rm rank} B\ge 1$. We are to construct a solution to match this upper bound. First, assume that $B$ is of full column rank, i.e.,${\rm rank} B=m\ge 1$, $m<n$. Moreover, assume that the first $m$ rows of $B$ are linearly independent, which can always be met by renumbering the states. We shall prove that, the perturbation pattern ${\cal A}=\left[ \begin{matrix}
{{0}_{m\times (n-m)}} & {{0}_{m\times m}}  \\
{{\bar I}_{n-m}} & {{0}_{(n-m)\times m}}  \\
\end{matrix} \right]$, ${\cal B}={{0}_{n\times m}}$, is feasible for to \ref{prob1}.
To see this, let $J^R=[m]$, $J^C=\{n+1,...,n+m\}$. For $i=1,...,p$, it holds that ${{\left[ {{\lambda }_{i}}I-A,B \right]}_{J^{R},J^{C}}}$, i.e., $B_{[m],[m]}$, is invertible. Let $\bar{J}^{R}=[n]\backslash J^{R}=\{m+1,...,n\}$ and $\bar{J}^{C}=[n-m]\subseteq [n+m]\backslash J^{C}$. Then, it turns out that $[{\cal A}, {\cal B}]_{\bar J^R,\bar J^C}$, which is exactly $I_{n-m}$, is of full generic rank. Hence, from Lemma 2, Condition (a) in Statement 2) of Theorem \ref{theorem 1} is satisfied.

Next, notice that ${\cal A}$ corresponds to that, a set of edges $\text{ }\!\!\{\!\!\text{ }({{z}_{1}},{{x}_{m+1}}),...,({{z}_{n-m}},{{x}_{n}})\}$ exists in the associated ACG ${\cal G}_{\rm auc}$. As argued in the proof of Theorem \ref{theorem 1}, a set of edges $\{({{x}_{1}},{{z}_{1}}),...,({{x}_{n}},{{z}_{n}})\}$ exist in ${\cal G}_{\rm auc}$. Recalling that the first $m$ rows of $B$ are linearly independent, there is at least one nonzero entry in the $i$th row of $B$, $i=1,...,m$. Moreover, every diagonal entry of ${\rm{adj}}(\lambda I - A)$ has degree $n-1$ for the variable $\lambda$, while the off-diagonal entry is either zero or has degree $n-2$, where ${\rm adj}(\cdot)$ takes the adjacency matrix. As a result, there is at least one nonzero entry in the $i$th row of ${{(\lambda I-A)}^{-1}}B=\frac{1}{\det (\lambda I-A)}\text{adj}(\lambda I-A)B$, for	$i=1,...,m$, since the unique polynomial with degree $n-1$ cannot be varnished by other terms in at least one entry of the $i$th row of $\text{adj}(\lambda I-A)B$. Therefore, vertices $z_1,...,z_m$ are input-reachable. If $m\ge n/2$, then $n-m\le m$. Due to the existence of $(z_i,x_{m+i})$ and $(x_{m+i},z_{m+i})$, $1\le i \le n-m\le m$, we obtain that $z_{m+1},...,z_{n}$ are input-reachable. Otherwise, $m< n/2$. Again, owing to the existence of $(z_1,x_{m+1})$ and $(x_{m+1}, z_{m+1})$, we get $z_{m+1}$ is input-reachable. Repeating such a process $n-2m-1$ times, we obtain successively that $z_{m+2},...,z_{m+n-2m}$ are input-reachable (see Fig. \ref{auc1} for illustration). This immediately leads to that, $z_{n-m+1},...,z_{n}$ are input-reachable. Hence, Condition (b) in Statement 3) of Theorem \ref{theorem 1} is satisfied. This indicates the above  $({\cal A}, {\cal B})$ is feasible for \ref{prob1}.

If ${\rm rank} B= r < m$, we can replace $B$ in the above argument with $r$ linearly independent columns  and neglect the other $m-r$ columns of $B$. It is easy to see that the remaining reasoning is still valid.


{\bf{Case II:}} If ${\rm rank} B=0$, add a nonzero entry to $B$. Then, according to the proof of {\bf Case I}, there exists a perturbation with  sparsity $n-1$ that makes the renewed system controllable. Hence, the upper bound for ${\rm rank} B=0$ is also valid.

We now prove the lower bound. Suppose that $[{\cal A},{\cal B}]$ is a feasible perturbation pattern for \ref{prob1}. Suppose further $A$ has $p$ distinct eigenvalues $\lambda_i|_{i=1}^p$.  Then, by Lemma \ref{lemma-int}, there are four sets $J_{i}^{R}\subseteq [n]$, $J_{i}^{C}\subseteq [n+m]$, $\bar{J}_{i}^{R}\subseteq [n]\backslash J_{i}^{R}$ and  $\bar{J}_{i}^{C}\subseteq [n\text{+}m]\backslash J_{i}^{C}$, such that $|J^R_i|=|J^C_i|$, $|\bar{j}_{i}^{C}|=|\bar{j}_{i}^{R}|= n-|j_{i}^{C}|$ and $[\lambda_iI -A, B]_{J_i^R,J_i^C}$ and ${{[{\cal A},{\cal B}]}_{\bar{j}_{i}^{R},\bar{j}_{i}^{C}}}$ are (generically) invertible, for each $i\in\{1,...,p\}$. Let $\bar{J}=\bigcup\nolimits_{i=1}^p{\bar{J}_{i}^{R}}$, and $k=\text{ }\!\!|\!\!\text{ }\bar{J}\text{ }\!\!|\!\!\text{ }$. Moreover, denote by $\tilde{\bar{J}}_{i}^{C}=\{n+m+1,...,n+m+|\bar{J}_{i}^{C}|\}$. Then, it turns out that both ${{[{{\lambda }_{i}}I-A,B,{{0}_{n\times n}}]}_{J_{i}^{R},J_{i}^{C}}}$  and ${{\left[ {{0}_{n\times n}},{{0}_{n\times m}},{{I}_{{\bar{J}}}} \right]}_{\bar{J}_{i}^{R}\tilde{\bar{J}}_{i}^{C}}}$  are invertible,  $i=1,...,p$. By Theorem \ref{theorem 1}, $\left[ {{0}_{n\times n}},{{0}_{n\times m}},{{I}_{{\bar{J}}}} \right]$ is feasible for \ref{prob1} associated with system $(A,[B,{{0}_{n\times n}}])$. That is, any feasible pattern $[{\cal A},{\cal B}]$ to \ref{prob1} corresponds to a pattern $\left[ {{0}_{n\times n}},{{0}_{n\times m}},{{I}_{{\bar{J}}}} \right]$ with sparsity $k\le {{\left\| [{{{{\cal A}}}_{o}},{{{{\cal B}}}_{o}}] \right\|}_{0}}$ that is feasible for Problem (\ref{lower bound}). Hence, the optimal value of Problem (\ref{lower bound}) is a lower bound of that of \ref{prob1}.  $\hfill\blacksquare$

Theorem \ref{bounds} gives an upper and a lower bound for \ref{prob1}. The lower bound means that, the ZNCD of a system is no less than the number of dedicated inputs (a dedicated input is an input that actuates only one state variable) that need to be added to the original system for achieving controllability. This builds a connection between ZNCD and the MCP.
The upper bound, $n-{\rm rank} B$, is the difference between the state dimension and the number of independent input vectors of $B$. Note that in \ref{prob1}, only the entries of the system matrices can be perturbed but the change of input number is not allowed.  Hence, this bound is nontrivial, especially when the number of inputs is limited. This bound indicates that the ratio between the number of perturbed entries for controllability and that of the total entries ($n^2+nm$) is upper-bounded by $O(\frac{1}{n+m})$, which approaches zero as $n$ increases. In what follows, We provide an example where the proposed bounds are tight.

\begin{example}[ZNDC of complete graphs] \label{fullone-1} For a complete graph with $n$ nodes (denoted by $K_n$), suppose the dynamics on it is characterized by the system matrices $A=1_{n\times n}$ and $B=1_{n\times 1}$. It can be calculated that ${\rm rank}({\cal C}(A,B))=1$, $\forall n>1$. Theorem \ref{bounds} predicts that $r_c(A,B)\le n-1$. Moreover, since ${\rm rank}([A,B])=1$, the optimal value to Problem (\ref{lower bound}) associated with $(A,B)$ is at least $n-1$. Hence, Theorem \ref{bounds} restricts $r_c(A,B)=n-1$.
\end{example}


\begin{remark}
	It is worth mentioning that \cite{Y_Zhang_2017} has shown that when $A$ and $B$ are both structured,  the minimal number of free entries that need to be added to achieve structural controllability equals the minimal number of dedicated inputs that need to be added for the same purpose. However, a similar relation does not exist between \ref{prob1} and Problem (\ref{lower bound}), although the optimal value of Problem (\ref{lower bound}) bounds that of \ref{prob1}; see the following example.
\end{remark}

\begin{example} \label{mid-example}
	Consider a system $(A,B)$ as
	$$A{\rm{ = }}\left[ {\begin{array}{*{20}{c}}
		{{\rm{ - }}1}&0&1&0&0&0\\
		0&{{\rm{ - }}1}&0&0&0&0\\
		0&0&3&0&0&0\\
		0&0&1&4&0&0\\
		0&0&0&0&4&0\\
		0&0&0&0&0&4
		\end{array}} \right], B=\left[
	\begin{array}{c}
	0 \\
	1 \\
	1 \\
	0 \\
	0 \\
	1 \\
	\end{array}
	\right].$$
	\ref{prob1} associated with $(A,B)$ has the optimal solution with cardinality $3$.
	This can be attained by verifying that, on the one hand, any perturbations with $2$ entries cannot make $(A,B)$ controllable. On the other hand, the $(1,1)$th, $(3,4)$th, $(5,6)$th entries of $A$ can be perturbed with the addition of $1$ to make the resulted system controllable. However, Problem (\ref{lower bound}) associated with $(A,B)$ has the optimal solution with cardinality $2$ (for example, two additional dedicated inputs actuating the $3$th and the $5$th states making the resulted system controllable).
\end{example}

	
	\section{Algorithms} \label{alg-sec}

	In this section, we provide two heuristic algorithms for computing ZNDC. The first one is a greedy algorithm, and the second one is based on the sequential convex relaxation. Both algorithms can deal with certain additional constraints imposed on the corresponding perturbations.
\subsection{Greedy Algorithm}
In the following, we develop a greedy algorithm for \ref{prob1} based on the structural controllability criterion in Theorem \ref{theorem 1}. This algorithm allows the perturbed entries to be chosen from a prescribed set. Let $ E_{zx}^{\rm ini}\subseteq  V_z\times V_x$ be the set corresponding to such entries. For an $E\subseteq E^{\rm ini}_{zx}$, define the pattern matrix  ${\cal M}_{E}\in \{0,*\}^{n\times (n+m)}$ as $[{\cal M}_{E}]_{ji}=*$ if and only if $(z_i,x_j)\in E$.

The key ingredient of our greedy algorithm is a suitable objective function, which serves as an `index' for controllability and indicates which local choice should be made greedily. To this end,  suppose $A$ has $p$ distinct eigenvalues $\lambda_i|_{i=1}^p$.   For each $i=1,...,p$, define the matroid ${\cal M}_i={\cal M}([I_n,\lambda_iI-A,B])$. Additionally, let the matroid ${\cal M}({\cal A}, {\cal B})={\cal M}(\left[[{\cal A}, {\cal B}]^\intercal, \bar I_{n+m} \right])$. Denote ${\bf M}_{n;m}=\{0,*\}^{n\times (n+m)}$. 
Define the function $g_1(\Delta A, \Delta B)$: $2^{{\bf M}_{n;m}}\to {\mathbb N}$:
\begin{equation}
g_1({\cal A}, {\cal B})=\sum \limits_{i=1}^p \rho({\cal M}_i \cap {\cal M}({\cal A}, {\cal B})).
\end{equation}
If $[{\cal A}, {\cal B}]={\cal M}_E$ for $E\subseteq E^{\rm ini}_{zx}$, we write $g_1({\cal M}_E)$ for $g_1({\cal A}, {\cal B})$ (the same below).

Moreover, the function $g_2({\cal A}, {\cal B})$: $2^{{\bf M}_{n;m}}\to {\mathbb N}$ is defined as the number of input-reachable vertices of $V_z$ in ${\cal G}_{\rm auc}$ associated with $({\cal A}, {\cal B})$. Since verifying whether a vertex is input-reachable can be done via the strongly connected component decomposition \cite{Geor1993Graph}, $g_2({\cal A}, {\cal B})$ can be computed in polynomial time. Additionally, the increase of $g_2({\cal A}, {\cal B})$ can be computed dynamically (or iteratively) at each stage when a new $e\in E^{\rm ini}_{zx}$ is added.

Based on the above, the objective function $g({\cal A}, {\cal B}):$ $2^{{\bf M}_{n;m}}\to {\mathbb R}$ for the greedy algorithm is defined as:
\begin{equation} \label{utility-fun}
g({\cal A}, {\cal B})=g_1({\cal A}, {\cal B})+\gamma \cdot g_2({\cal A}, {\cal B}),
\end{equation}
where $\gamma>0$ is a given constant that weights the relative importance between $g_1({\cal A},{\cal B})$ and $g_1({\cal A},{\cal B})$. In particular, if $\gamma\gg 1$, then the greedy algorithm tends to select $e\in E_{zx}^{\rm ini}$ with bigger increase in $g_2({\cal A}, {\cal B})$ in the early stages; if $\gamma \ll 1$, the contrary. Typically, $\gamma =1$. It is easy to see that, for any $\gamma>0$, $g({\cal A}, {\cal B})$ is non-decreasing over $2^{{\bf M}_{n;m}}$. Moreover, $g({\cal A}, {\cal B})=pn + \gamma (n+m)$ means the system in (\ref{state space perturbation}) is structurally controllable.

	\begin{algorithm} 
		{{{{
						\caption{: A greedy algorithm for \ref{prob1} with a prescribed set of perturbable entries ${\cal M}_{E^{\rm ini}_{zx}}$} 
						\label{alg1} 
						\begin{algorithmic}[1] 
							\STATE Calculate the eigenvalues $\lambda_i|_{i=1}^p$ of $A$, and construct the ACG ${\cal G}_{\rm auc}(0_{n\times (n+m)})$.
							\STATE Initialize $E\leftarrow \emptyset$
							\WHILE{$g({\cal M}_E)  <  pn+\gamma (n+m) $}
							\STATE $e\leftarrow $ $e' \in \arg \max \nolimits_{a\in E^{\rm ini}_{zx}\backslash E} \ g({\cal M}_{E\cup \{a\}}) - g({\cal M}_{E})$.
							\STATE $E\leftarrow E\cup \{e\}$;
							\ENDWHILE
							\STATE Return ${\cal M}_E$.
				\end{algorithmic}}}
		}}
	\end{algorithm}
	
	It is worth noting that Algorithm \ref{alg1} is heuristic without optimality guarantee. Moreover, due to the non-submodular or non-supermodular of $g({\cal A}, {\cal B})$,  there seem no nontrivial approximation bounds. Nevertheless, Algorithm \ref{alg1} performances fairly well in practice as shown in Section \ref{simu-sec}.

\subsection{Algorithm based on sequential convex relaxation}
In this section, we provide a heuristic algorithm for \ref{prob2} based on the weighted $l_1$-norm relaxation and sequential convex optimization. More precisely, we first relax \ref{prob2} by approximating the $||\cdot ||_0$ with a weighted $l_1$-norm. Then, we rewrite the controllability constraint as a rank constraint involving the controllability Gramian and the Lyapunov function. The rank constraint is subsequently reformulated as minimizing the difference between two Ky Fan norms (inspired by \cite{doelman2016sequential,hu2012fast}), leading to a difference-of-convex-function problem (DCP) combined with the weighted $l_1$-norm. The DCP is then solved by the standard concave-convex procedure, which is assured to converge. 

To formulate the DCP of \ref{prob2}, the following intermediate result is first presented.

\begin{lemma}\cite{Geor1993Graph} \label{circle-theo} Let $\rho(X)$ be the spectrum of $X=[x_{ij}]$, i.e., the maximum magnitude of eigenvalues of $X$. Then
	$\rho(X)\le \min\{ \max \limits_{i} \sum \limits_{j=1}^n |x_{ij}|, \max  \limits_{j} \sum \limits_{i=1}^n |x_{ij}|\}$.
\end{lemma}

At first, since the $||\cdot||_0$ in \ref{prob2} is non-convex and non-smooth, an usual heuristic is to use $||\cdot ||_1$ ($l_1$-norm) to approximate $||\cdot||_0$. However, if $\theta$ is feasible for \ref{prob2}, then for any closed interval $[\kappa_1,\kappa_2]$, there are infinitely many $\kappa\in [\kappa_1,\kappa_2]$ such that $\kappa \theta$ is also feasible due to the genericity of controllability. Instead of using $||\cdot||_1$, inspired by \cite{sriperumbudur2011majorization}, we adopt a weighted $l_1$-norm for $\theta\in {\mathbb R}^{l}$  as
\begin{equation}
||\theta||_{[\tau]} \doteq \sum \limits_{i=1}^l \frac{\log(1+|\theta_i|/\tau)}{\log(1+1/\tau)}, \ \tau>0.
\end{equation}
It is easy to verify that
\begin{equation}
||\theta||_0 = \lim \limits_{\tau\rightarrow 0_+} \sum \limits_{i=1}^l \frac{\log(1+|\theta_i|/\tau)}{\log(1+1/\tau)}.
\end{equation}
Indeed, for $\theta_i\ne 0$, $\lim \limits_{\tau\rightarrow 0_+} \frac{\log(1+|\theta_i|/\tau)}{\log(1+1/\tau)}= \lim \limits_{\tau\rightarrow 0_+}\frac{\log(|\theta_i|\tau^{-1})}{\log\tau^{-1}} =  \lim \limits_{\tau\rightarrow 0_+}\frac{\log |\theta_i|+\log \tau^{-1}}{\log \tau^{-1}} =1$. It has been pointed out in \cite{sriperumbudur2011majorization} that, compared to $||\theta ||_1$, $||\theta||_{[\tau]}$ is a tighter approximation to $||\theta||_0$ for any $\tau>0$. Therefore, sparser solutions are expected for \ref{prob2} obtained by replacing $||\theta||_0$ with $||\theta||_{[\tau]}$, compared to replacing $||\theta||_0$ with $||\theta||_{1}$.

To formulate the controllability constraint, we adopt the Lyapunov equation. From \citep[Lemma 3.18]{zhou1996robust}, for $(A,B)$ with $A$ stable (i.e., all eigenvalues of $A$ have negative real parts), $(A,B)$ is controllable, if and only if the solution to the Lyapunov equation
$$AW+WA^\intercal+BB^{\intercal}=0$$
is positive definite (in fact, the solution $W$ is the controllability Gramian of $(A,B)$). With this idea and the weighted $l_1$-norm, we consider a relaxation of \ref{prob2} as

\begin{align}  \label{prob4}
\mathop {\min }\limits_{\theta\in {\mathbb R}^l, W\in {\mathbb S}^n} &{\kern 1pt} {\kern 1pt} {\kern 1pt} {\kern 1pt} {\kern 1pt} \sum \limits_{i=1}^l \frac{\log(1+ \theta_i/\tau)}{\log(1+1/\tau)}\\
{\rm s.t.} \ \ &(A(\theta)-\mu I)W+W(A^\intercal(\theta)-\mu I)+B(\theta)B^{\intercal}(\theta)=0 \label{const1} \\
& W\succeq \varepsilon I \label{const2} \\
& \theta \le 1_{l\times 1} \label{const3} \\
& \theta  \ge 0_{l\times 1} \label{const4} \\
& \sum \nolimits_{i=1}^l \theta_i \ge 1-\eta \label{const5}
\end{align}
where $0<\tau \ll 1$, $\mu> \max \limits_{i} \{\sum \nolimits_{j=1}^n (|A_{ij}|+\sum \nolimits_{k=1}^l |[A_{k}]_{ij}|)\}$, $\varepsilon >0$ but is arbitrarily close to $0$, and $0<\eta<1$. Here, (\ref{const5}) is introduced so that the optimal $\theta$ will not approach to zero ($\eta$ can be close to zero), while $\mu$ is introduced in (\ref{const1}) so that $A(\theta)-\mu I$ is stable for $\theta$ satisfying (\ref{const3})-(\ref{const5}). In practice, $\mu$ can be much bigger than $\max \limits_{i} \{\sum \nolimits_{j=1}^n (|A_{ij}|+ \sum \nolimits_{k=1}^l |[A_{k}]_{ij}|)\}$.  The introduction of constraints (\ref{const3})-(\ref{const5}) is justified in the following lemma.

\begin{lemma}  If \ref{prob2} has a feasible solution $\theta$ with $||\theta||_0\le l^*$, $l^*\in {\mathbb N}$, then Problem (\ref{prob4}) also has a feasible solution $\theta'$ satisfying $||\theta'||_0\le l^*$, provided that $\varepsilon$ is arbitrarily close to zero and $\mu> \max \limits_{i} \{\sum \nolimits_{j=1}^n (|A_{ij}|+\sum \nolimits_{k=1}^l |[A_{k}]_{ij}|)\}$.
\end{lemma}
\begin{proof} According to Lemma \ref{circle-theo}, for $\theta$ subject to (\ref{const3})-(\ref{const5}),  $\rho(A(\theta))\le \max \nolimits_{i}\{ \sum \nolimits_{j=1}^n (|A_{ij}|+\sum \nolimits_{k=1}^l \theta_k |[A_k]_{ij}|)\} \le \max \nolimits_{i} \{\sum \nolimits_{j=1}^n (|A_{ij}|+\sum \nolimits_{k=1}^l |[A_k]_{ij}|)\}$, where the second inequality is due to (\ref{const3}). Therefore, all eigenvalues of $A(\theta)-\mu I$ have negative real parts, i.e., $A(\theta)-\mu I$ is stable. On the other hand, since controllability is a generic property and the set for $\theta$ subject to (\ref{const3})-(\ref{const5}) is dense, if $(A(\theta), B(\theta))$ is structurally controllable, there exists a $\theta'=\kappa \theta$ for some $\kappa \in {\mathbb R}$ so that $\theta'$ satisfies (\ref{const3})-(\ref{const5}) and $(A(\theta'), B(\theta'))$ is controllable, meanwhile $||\theta'||_0=||\theta||_0$. From the PBH test, the controllability of $(A(\theta'), B(\theta'))$ implies the controllability of $(A(\theta')-\mu I, B(\theta'))$. Noting $A(\theta')-\mu I$ is stable, the Lyapunov equation (\ref{const1}) has a positive definite solution $W$. This finishes the proof.
\end{proof}

Notice that (\ref{const1}) is bilinear in the variables $\theta$ and $W$. To handle this, inspired by \cite{doelman2016sequential}, we introduce an equivalent rank constraint. To this end, define matrices
$$M\doteq [A(\theta)-\mu I, W, B(\theta)], N\doteq [W, A(\theta)-\mu I, B(\theta)]^\intercal.$$Then, (\ref{const1}) is equivalent to $MN=0$. Introduce the matrix variable $Z\in {\mathbb R}^{(3n+m)\times (3n+m)}$ as
\begin{equation}
Z= \left[
\begin{array}{cc}
0 & M \\
N & I_{2n+m} \\
\end{array}
\right].
\end{equation}
According to the Schur complement,
${\rm rank}(Z)={\rm rank}(I_{2n+m})+ {\rm rank}(MN)$. 	Hence, $MN=0$, if and only if ${\rm rank}(Z)=2n+m$.

The rank constraint ${\rm rank}(Z)=2n+m$ is still hard to handle. Inspired by \cite{hu2012fast},  we can replace the rank constraint with the truncated nuclear norm as follows:
\begin{equation} \label{truncate-norm}
||Z||_{*}-||Z||_{F_{2n+m}}=0,
\end{equation}
where $||Z||_*$ denotes the nuclear norm of $Z$, i.e., the sum of its nonzero singular values, and $||Z||_{F_{2n+m}}$ denotes the Ky Fan $(2n+m)$-norm, namely, the sum of the largest $2n+m$ singular values of $M$ (the Ky Fan $r$-norm is defined similarly for any $r\in {\mathbb N}_+$).  Keeping in mind of (\ref{truncate-norm}), we can relax Problem (\ref{prob4}) as the following difference of convex problem (DCP) \cite{yuille2003concave,lanckriet2009convergence}
\begin{align}   \label{prob5}
\mathop {\min }\limits_{\theta, W, Z} & \ F(\theta, W, Z)\doteq \sum \limits_{i=1}^l \frac{\log(1+ \theta_i/\tau)}{\log(1+1/\tau)}+ \gamma (||Z||_{*}-||Z||_{F_{2n+m}})\\
{\rm s.t.}\ & Z=\left[
\begin{array}{cc}
0 & M \\
N & I_{2n+m} \\
\end{array}
\right] \label{const-add}  \\
& (\ref{const2})-(\ref{const5}) \label{const6}
\end{align}
where $\gamma>0$ is a regularization parameter. The DCP of Problem (\ref{prob5}) comes from the fact that the objective can be rearranged as $F(\theta, W, Z)=\gamma||Z||_{*}+ \sum \nolimits_{i=1}^l \frac{\log(1+ \theta_i/\tau)}{\log(1+1/\tau)}-\gamma ||Z||_{F_{2n+m}}$, where the first item is convex and the last two items are concave in the minimization variables $(W,\theta, Z)$ (the nuclear norm and Ky Fan $r$-norm for any $r\in {\mathbb N}$ are convex in the matrix set \cite{hu2012fast}).  Note also that $Z$ is affine in $\theta$ and $W$ from (\ref{const-add}). Therefore, the constraints (\ref{const-add})-(\ref{const6}) define a convex set for the minimization variables $(W,\theta, Z)$. It can be seen that, for $\tau\rightarrow 0$ and $\varepsilon\rightarrow 0$, the optimal solution $\theta^*$ to ${\cal P}_2$ also corresponds to an optimal $(W^{*'}, \theta^{*'}, Z^{*'})$ that minimizes Problem (\ref{prob5}). But in practice, the parameter $\tau$ cannot be arbitrarily close to zero due to the limitation of computer precisions. 

\begin{remark}We may also drop the decision variable $Z$ from Problem (\ref{prob5}) without breaking its property of being a DCP, since $Z$ is affine in $(\theta, W)$. However, regarding $Z$ as a decision variable can benefit us in deriving the first-order derivative of the concave part in the objective $F(\theta, W, Z)$.
\end{remark}

The standard method for DCPs is the well-established concave-convex procedure (CCCP) \cite{yuille2003concave}, which returns locally optimal solutions via solving a sequence of convex programs with convergence guarantees. The key of the CCCP is to determine the first-order derivative of the concave part in the objective to linearize it around a solution $(W^{(k)}, \theta^{(k)}, Z^{(k)})$ obtained in the current ($k$th) iteration. Towards this end, we establish
\begin{equation} \label{linear-appro}
\begin{array}{c}	\sum \limits_{i=1}^l \frac{\log(1+ \theta_i/\tau)}{\log(1+1/\tau)}\thickapprox \sum \limits_{i=1}^l \left\{\frac{
	\log(1+ \theta^{(k)}_i/\tau)}{\log(1+1/\tau)}+ \right. \\
\left.	\frac{\tau^{-1}}{(\log(1+\tau^{{-1}}))(1+\theta^{(k)}_i\tau^{{-1}})}(\theta_i-\theta^{[k]}_i) \right\}\end{array},
\end{equation} 
Additionally, let the singular value decomposition of $Z^{(k)}$ be
\begin{equation}\label{svd} Z^{(k)}=[U_1^{(k)}, U_2^{(k)}]\left[
\begin{array}{cc}
\Lambda_1 &  \\
& \Lambda_2 \\
\end{array}
\right]\left[
\begin{array}{c}
V^{(k),\intercal}_1 \\
V^{(k),\intercal}_2 \\
\end{array}
\right]
,\end{equation} where $U_1^{(k)}$ and $V_1^{(k)}$ are respectively the left and right singular vectors associated with the largest $2n+m$ singular values of $Z^{(k)}$. From \citep[Theorem 3.4]{qi1996extreme}, $U_1^{(k)}V_1^{(k),\intercal}\in  \partial ||Z||_{F_{2n+m}}$ ($\partial$ denotes the subdifferential). Therefore, $-||Z||_{F_{2n+m}}$ can be linearized at $Z^{(k)}$ as \cite{watson1992characterization,qi1996extreme}
\begin{equation}\label{Fr-norm}-||Z||_{F_{2n+m}}\thickapprox -||Z^{(k)}||_{F_{2n+m}}-{\rm tr}(U_1^{(k),\intercal}(Z-Z^{(k)})V_1^{(k)}),\end{equation} where ${\rm tr}(\cdot)$ takes the trace.


Based on the above expressions, the convex program at the $(k+1)$th iteration is formulated as
	
\begin{equation} \label{sub-convex}
\begin{aligned}
\mathop {\min }\limits_{\theta, W, Z} &{\kern 1pt} {\kern 1pt} {\kern 1pt} {\kern 1pt} {\kern 1pt} \gamma||Z||_{*}+ \\
& \sum\limits_{i=1}^l \frac{\tau^{-1}}{(\log(1+\tau^{{-1}}))(1+\theta^{(k)}_i\tau^{{-1}})}\theta_i-\gamma {\rm tr}(U_1^{(k),\intercal}ZV_1^{(k)})\\
{\rm s.t.} \ \ \  &  (\ref{const-add}), (\ref{const6})
\end{aligned}
\end{equation}
For ease of explication, we collect the procedure for approximating \ref{prob2} as Algorithm \ref{alg2}. The convergence of Algorithm \ref{alg2} is stated in the following theorem.

\begin{algorithm} 
	{{{
				\caption{: A convex-relaxation algorithm for \ref{prob2}} 
				\label{alg2} 
				\begin{algorithmic}[1] 
					\STATE {Initialize $k=0$ and $\theta^{(0)}, W^{(\theta)}$ and $Z^{(0)}$ by solving the convex problem
						\begin{equation}  \label{initial-prob}
						\begin{aligned}
						& \mathop {\min }\limits_{\theta,W,Z} \ 1 \\
						{\rm s.t.} \ & \ (\ref{const2})-(\ref{const5}), (\ref{const-add})
						\end{aligned}
						\end{equation}
					}
					\WHILE {$||\theta^{(k)}-\theta^{(k-1)}||>\xi$ ($\xi>0$ is the convergence threshold)} 
					\STATE Obtain $U_1^{(k)}$ and $V_1^{(k)}$ according to (\ref{svd});
					\STATE Solve the convex program (\ref{sub-convex}) to obtain $(\theta^{(k+1)},W^{(k+1)}, Z^{(k+1)})$;
					\STATE $k+1 \leftarrow k$;
					\ENDWHILE
					\STATE Return $\theta^{(k)}$ when converging.
		\end{algorithmic}}}
	}
\end{algorithm}

\begin{theorem}The sequence $\{ \theta^{(k)}, W^{(k)}, Z^{(k)}\}$ generated by Algorithm \ref{alg2} satisfies:
	
	(i) $F(\theta^{(k+1)}, W^{(k+1)}, Z^{(k+1)})\le F(\theta^{(k)}, W^{(k)}, Z^{(k)})$;
	
	(ii) $\lim \limits_{k\rightarrow \infty} (F(\theta^{(k+1)}, W^{(k+1)}, Z^{(k+1)})-F(\theta^{(k)}, W^{(k)}, Z^{(k)}))=0$ and $\lim \limits_{k\rightarrow \infty} (\theta^{(k+1)}-\theta^{(k)})=0$.
	
	Moreover, Algorithm \ref{alg2} is guaranteed to converge to a stationary point of Problem (\ref{prob5}).
\end{theorem}

\begin{proof}
	Note Algorithm \ref{alg2} is a CCCP for Problem (\ref{prob5}) that is a DCP. Hence, property (i) comes immediately from the non-increasing property of the CCCP \citep[Theorem 2]{yuille2003concave}. For property (ii) and the convergence statement, note the `$\approx$' in (\ref{linear-appro}) and (\ref{Fr-norm}) can be replaced by $\le$, for any $Z\in {\mathbb R}^{(3n+m)\times (3n+m)}$, $\theta \in {\mathbb R}^{l}$, due to the concavity. Hence, Algorithm \ref{alg2} is equivalent to the majorization-minimization (MM) algorithm (see \citep[Section 2]{lanckriet2009convergence}) by replacing the objective function in Problem (\ref{prob5}) with the corresponding approximate function, which is obtained as the sum of $||Z||_{*}$ and the right-hand side of (\ref{linear-appro}) and (\ref{Fr-norm}) with the respective regularization parameters. Then, it is easy to verify that this approximate function satisfies the sufficient conditions in \citep[Theorem \ref{theorem 1} \& Corollary 1]{razaviyayn2013unified}, and according to that\footnote{\citep[Theorem \ref{theorem 1}]{razaviyayn2013unified} states the limit point returned by the MM algorithm is the stationary point solution of the original problem if the approximate function satisfies: 1) it is continuous, 2) it is a tight upper bound of the original objective function, and 3) it has the same first-order directional derivative of the original objective function at the point where the upper bound is tight.}, the MM algorithm is guaranteed to converge to a stationary point of Problem (\ref{prob5}). This proves property (ii) and the convergence statement.
\end{proof}



Recall a {\emph{stationary point}} refers to a point which satisfies the corresponding KKT conditions of the optimization problem (necessary for optimality) \cite{lanckriet2009convergence}. It is noted that there is no guarantee that Algorithm \ref{alg2} returns an {\emph{optimal solution}} to Problem (\ref{prob5}), not to mention the optimal solution to the {\emph{original}} \ref{prob2}.

\begin{remark} For $Z\in {\mathbb R}^{n_1\times n_2}$, the nuclear norm $||Z||_*$ can be expressed as \citep[Appendix A]{nguyen2019low} $||Z||_*=\min \{\frac{1}{2}{\rm tr}(W_1)+\frac{1}{2}{\rm tr}(W_2): \left[
	\begin{array}{cc}
	W_1 & Z \\
	Z^\intercal & W_2 \\
	\end{array}
	\right]\succeq 0\}
	$. Hence, the convex program (\ref{sub-convex}) can be converted into the following semi-definite program:
	\begin{equation}  \label{sub-SDP}
	\begin{aligned}
	\mathop {\min }\limits_{\theta, W, Z, W_1, W_2} &{\kern 1pt} {\kern 1pt} {\kern 1pt} {\kern 1pt} {\kern 1pt} \frac{1}{2}\gamma{\rm tr}(W_1+W_2)+ \\
	& \sum\limits_{i=1}^l \frac{\tau^{-1}}{\log(1+\tau^{{-1}})(1+\theta^{(k)}_i\tau^{{-1}})}\theta_i-\gamma {\rm tr}(U_1^{(k),\intercal}ZV_1^{(k)})\\
	{\rm s.t.} \ & \left[
	\begin{array}{cc}
	W_1 & Z \\
	Z^\intercal & W_2 \\
	\end{array}
	\right]\succeq 0, \ {\rm and} \  (\ref{const-add}), (\ref{const6})
	\end{aligned}
	\end{equation}
\end{remark}


\begin{remark}To accelerate the convergence of Algorithm \ref{alg2}, the objective of Problem (\ref{initial-prob}) in the initialization step can be replaced by $||W||_*$.
\end{remark}

\begin{remark}The simulations in the next section show that Algorithm \ref{alg2} could return sparse solutions for $\theta^{(k)}$. In practice, if $\theta^{(k)}_i\ll (1-\eta)/l$, it could be thought $\theta^{(k)}_i=0$.
\end{remark}



\section{Typical Examples and Simulations} \label{simu-sec}

In this section, we validate the effectiveness of the proposed algorithms in computing ZNDCs of several typical uncontrollable networks arising in multi-agent systems. Specially, for a network (graph) ${\cal G}=(V,E)$, its dynamics is characterized by (\ref{StateSpace}), with the corresponding state matrix $A$ satisfying $A_{ij}\ne 0$ only if the edge $(j,i)\in E$ ($i\ne j$, $i,j\in V$), corresponding to that, each node of ${\cal G}$ is of first-order dynamics. Each column of $B$ in (\ref{StateSpace}) corresponds to an input node.

\begin{example}[Example \ref{fullone-1} cont.] \label{fullone} Consider the complete graph $K_n$ in Example \ref{fullone-1}. From that example, we know the ZNCD of $K_n$ is $n-1$. The greedy algorithm (Algorithm \ref{alg1}) always returns the correct solution $r_c=n-1$ (in which the regularization parameter $\gamma=1$). An optimal perturbation structure is given as
	$$\left[
	\begin{array}{cc}
	\bar I_{n-1} & 0_{(n-1)\times 2} \\
	0_{1\times (n-1)} & 0_{1\times 2} \\
	\end{array}
	\right].
	$$
	For $n=6$, other possible optimal perturbation structures returned by Algorithm \ref{alg1} are illustrated in Fig. \ref{example-complete-graph}.
	
	
	\begin{figure} \centering
		\includegraphics[width=3.5in]{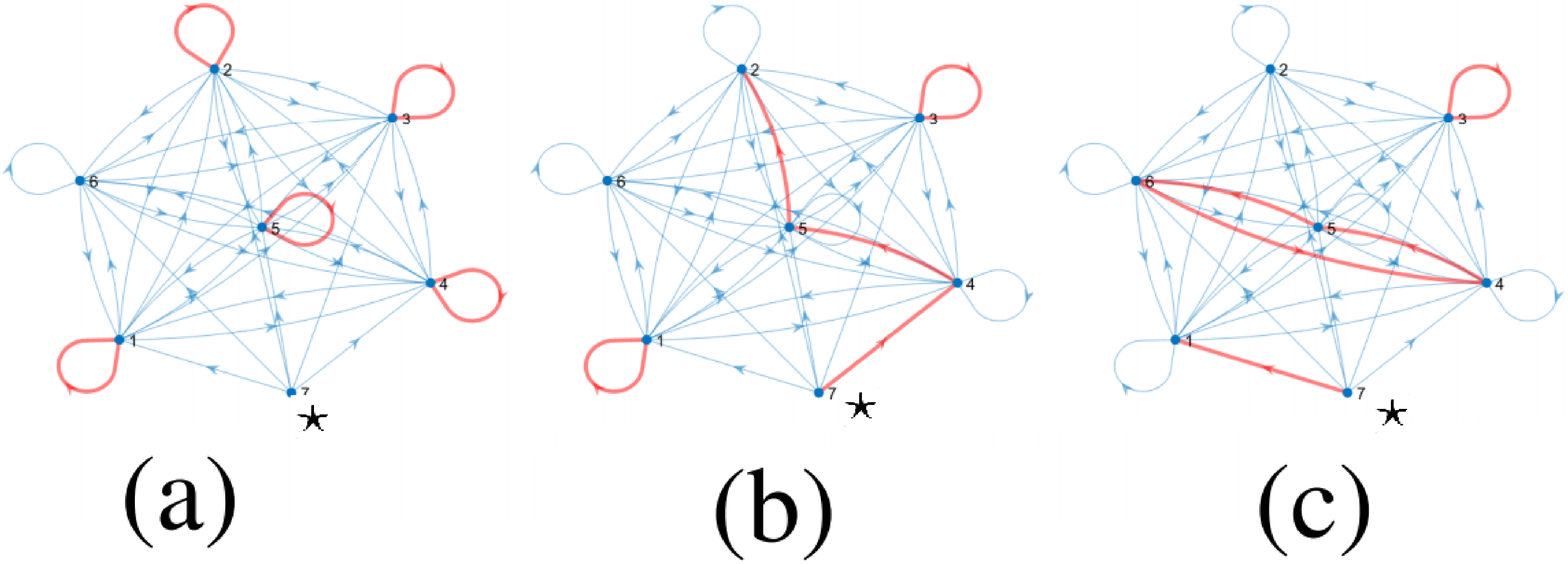}\\
		\caption{ZNDC of $K_6$ in Example \ref{fullone}. The node labeled $\star$ is the input node, and the rest are the state nodes. The same below.  Red edges indicate that those edges need to be perturbed to make the original system controllable in a solution returned by Algorithm \ref{alg1}. }
		\label{example-complete-graph}
	\end{figure}
\end{example}

\begin{example}[Generating structurally controllable networks] \label{construct} We show the application of Algorithm \ref{alg1} in constructing `minimal' structurally controllable networks. At first, consider $A={ 0}_{n\times n}$, $B={ 0}_{n\times 1}$, and all entries of $[A,B]$ can be perturbed. In this case, \ref{prob1} reduces to determining a structure with the minimum number of edges to ensure structural controllability. From the classical structural controllability theory \cite{Lin_1974}, we know such a structure is a path from the input node that walks through all the state nodes (i.e., a stem containing all the state nodes). We present the structurally controllable structures returned by Algorithm \ref{alg1} ($\gamma =1$) in Fig. \ref{example-path} (a)-(c), which are indeed the so-called paths. This validates the effectiveness of Algorithm \ref{alg1}. Next, suppose $A={0}_{n\times n}$, $B={ 0}_{n\times 1}$, and the perturbable entries belong to a prior set $E_{zx}^{\rm ini}$. We randomly generate $E_{zx}^{\rm ini}$, and the resulting structurally controllable networks returned by Algorithm \ref{alg1} ($\gamma=1$) are presented in Fig. \ref{example-path}(d)-(f). Through exhausting search (note this problem is NP-hard \cite{zhang2019minimal,Y_Zhang_2017}), it can be found that those networks are indeed the optimal ones.
	
	\begin{figure} \centering
		\subfigure { \label{examp1-a} 
			\includegraphics[width=0.54\columnwidth]{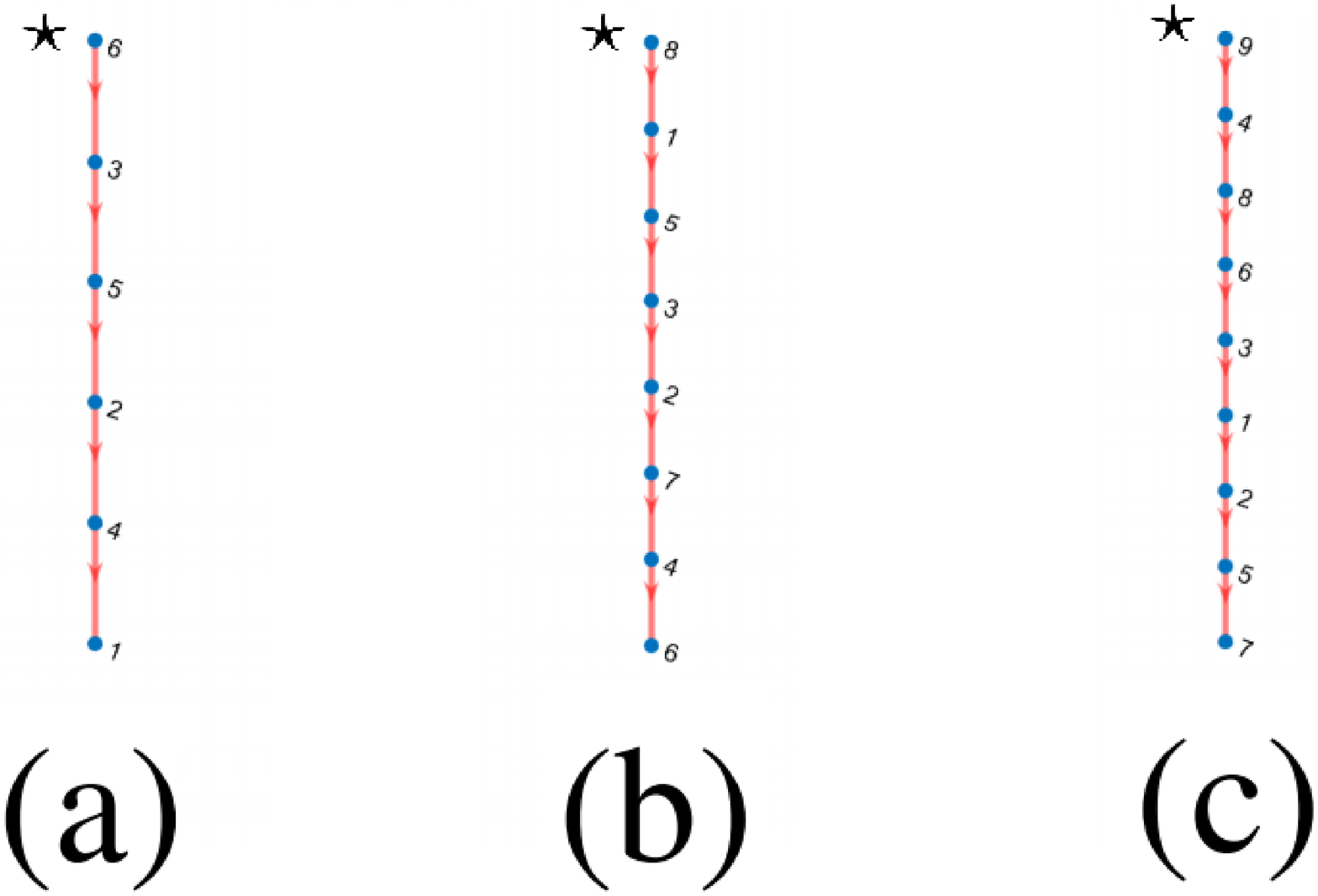}
		}
		\subfigure { \label{examp1-hl}
			\includegraphics[width=0.95\columnwidth]{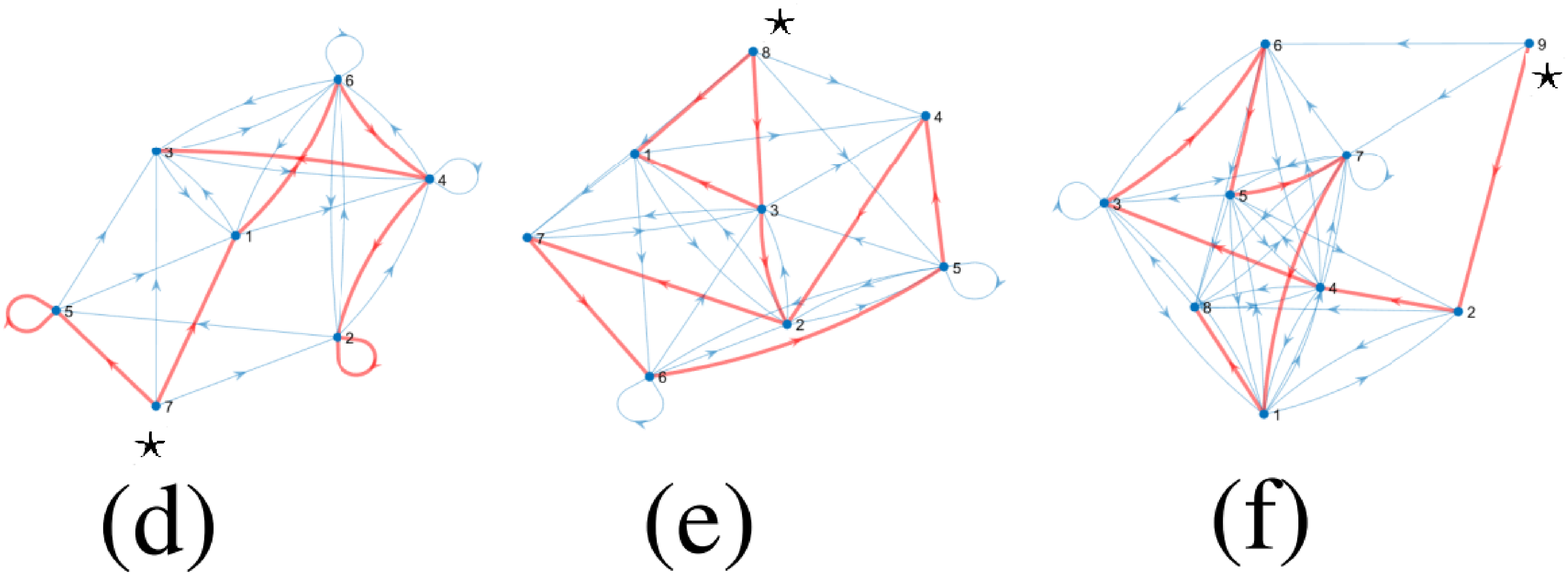}
		}
		\caption{Structurally controllable networks returned by Algorithm \ref{alg1} in Example \ref{construct}. Blue edges represent the edges that can be chosen (perturbed), while red ones represent the chosen edges via Algorithm \ref{alg1}. } 
		\label{example-path}
	\end{figure}
	
\end{example}

\begin{figure}
\centering
\includegraphics[width=3.5in]{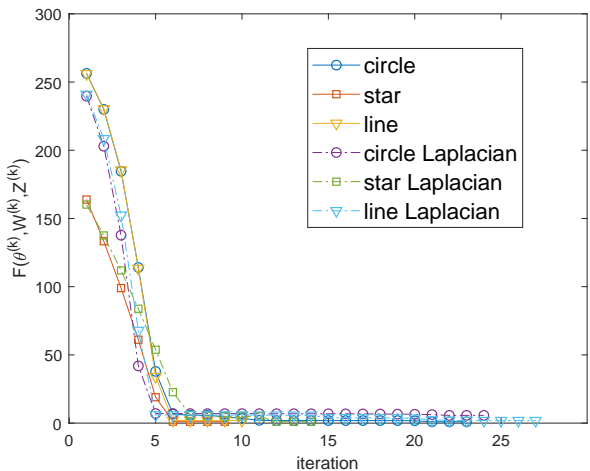}\\
\caption{Evolution of $F(\theta^{(k)}, W^{(k)}, Z^{(k)}))$ during the iteration.}\label{convergence}
\end{figure}

	\begin{example}\label{three-type} 
		%
		%
		%
	We continue to consider three typical networks: the line network ${\cal G}_{\rm line}$, star network ${\cal G}_{\rm star}$, and the circle network ${\cal G}_{\rm circ}$, shown respectively in Fig \ref{example-three-typical}(a)-(c), (d)-(f), and (g)-(i), all with $7$ state nodes and one input node. In those networks, all the edges have weight $1$, except the self-loops which have weight $-1$; besides, only the node indexed by $1$ has an external input. Denote the corresponding state matrix and input matrix by $A({\cal G}_\star)$ and $B_\star$, respectively, $\star={\rm star}, {\rm line}, {\rm circ}$. It can be calculated that ${\cal C}(A({\cal G}_{\rm star}),B_{\rm star})=2$, ${\cal C}(A({\cal G}_{\rm line}),B_{\rm line})=4$, ${\cal C}(A({\cal G}_{\rm circ}),B_{\rm circ})=4$, indicating that all those networks are uncontrollable. The uncontrollability of those networks is due to the symmetries in network topologies \cite{chapman2014symmetry}. Applying Algorithm \ref{alg1} to those networks by setting $\gamma=1$, we consider two scenarios: i) every entry of $[A({\cal G}_{\star}), B_\star]$ can be perturbed, and ii) only the existing edges in ${\cal G}_{\star}$ can the perturbed (i.e., $E_{zx}^{\rm ini}=E({\cal G}_{\star})$ with $\star={\rm star}, {\rm line}, {\rm circ}$). The outputs of Algorithm \ref{alg1} are given in Fig \ref{example-three-typical}(a)-(b), (d)-(e), and (g)-(h), respectively. From them, we observe that symmetries (in the edge weights) are broken in the controllable topologies. By the lower bound in Theorem \ref{bounds}, it can be validated that those solutions are optimal.


Next, suppose all edges in ${\cal G}_{\star}$ except the self-loops are undirected; that is, ${{[A({\cal G}_{\star})]}_{ij}}={{[A({\cal G}_{\star})]}_{ji}}$ for $i,j=1,...,n$, $\star = {\rm star}, {\rm line}, {\rm circ}$.  In this scenario, $[A({\cal G}_{\star}),B_\star]$ can be parameterized by the weights of the undirected edges of ${\cal G}_{\star}$ in the form of (\ref{affine}) (but not (\ref{Linear Parameterization}), since the coefficient matrices may have a rank of $2$). We adopt Algorithm \ref{alg2} for this scenario. The parameters of Algorithm \ref{alg2} are chosen as in Table \ref{parameters}. For all three networks, we plot the evolution of $F(\theta^{(k)}, W^{(k)}, Z^{(k)})$ versus $k$ in Fig. \ref{convergence} (labelled `{\rm star}, {\rm line}, {\rm circ}'). From this figure, the non-increasing of $F(\theta^{(k)}, W^{(k)}, Z^{(k)}))$, as well as the convergence of Algorithm \ref{alg2}, is validated. The solutions for the vector of edge weights when converging are {\emph{sparse}}, which are respectively

$\theta=[ 0.45, 0.00,    0.00,   0.45,    0.00,   0.00]^\intercal$ for ${\cal G}_{\rm line}$

$\theta=[0.00 ,   0.00,   0.00,   0.00 ,    0.00,   0.90]^\intercal$ for ${\cal G}_{\rm star}$

$\theta=[ 0.00,    0.00,   0.00,    0.00,   0.90,    0.00,    0.00]^\intercal$ for ${\cal G}_{\rm circ}$,
where the $i$th entry of $\theta$ corresponds to the $i$th nonzero entry in the upper triangular part of $A({\cal G}_{\star})$ (in the columnwise). Accordingly, the corresponding perturbed edges are given in Fig. \ref{example-three-typical}(c), (f), and (i), respectively. Note that Algorithm \ref{alg2} is not guaranteed to return a feasible solution to \ref{prob2} (the solution in Fig. \ref{example-three-typical}(f) is not feasible), since it is based on the relaxation Problem (\ref{prob5}) of \ref{prob1}. Compared to Algorithm \ref{alg1}, Algorithm \ref{alg2} is applicable to a larger class of linear parameterizations.


\begin{table}[htbp]
	\centering  
	\caption{Parameter setting for Algorithm \ref{alg2}}  
	\label{parameters}  
	\begin{tabular}{|c|c|c|c|c|}
		\hline  
		& & & &\\[-6pt]  
		$\tau$ & $\gamma$ & $\eta$ & $\varepsilon$ & $\xi$\\  
		\hline
		& & & &\\[-6pt]  
		$10^{-5}$& $40$ &$0.1$ &$10^{-5}$&$10^{-5}$\\
		\hline
	\end{tabular}
\end{table}

\begin{figure} \centering
	\subfigure { \label{examp1-a} 
		\includegraphics[width=0.93\columnwidth]{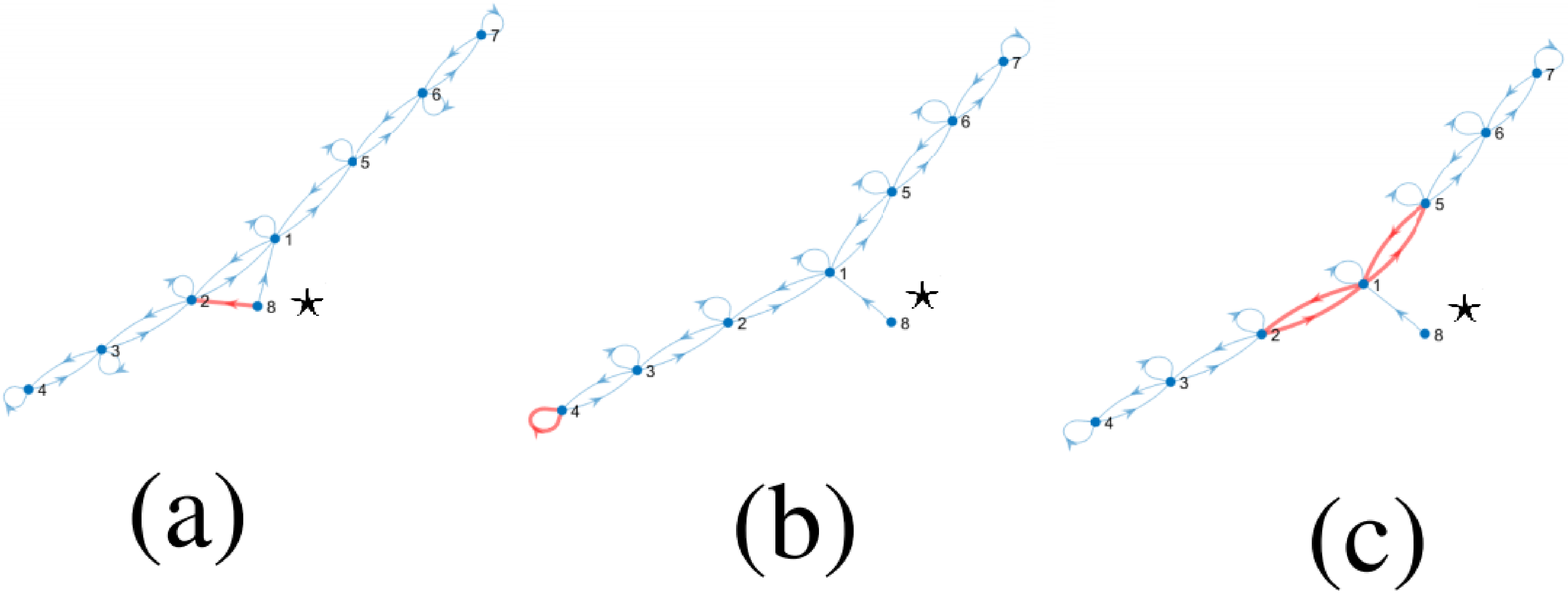}
	}
	\subfigure { \label{examp1-hl}
		\includegraphics[width=0.93\columnwidth]{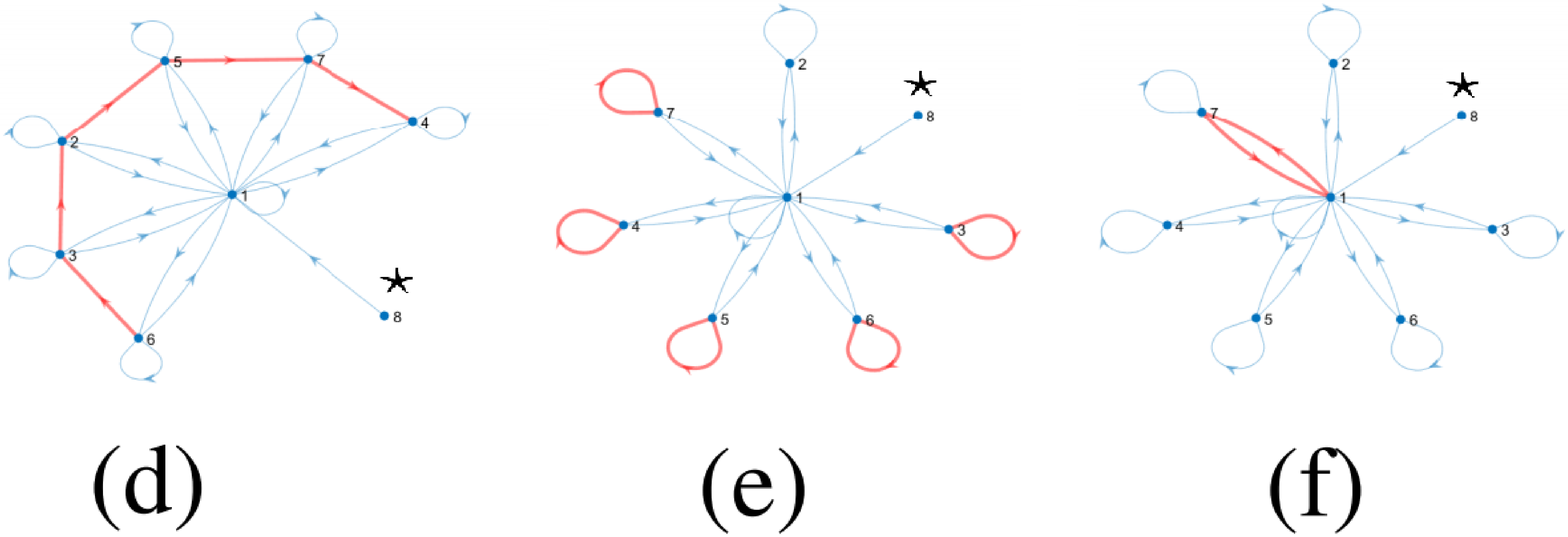}
	}
	\subfigure { \label{examp1-b}
		\includegraphics[width=0.93\columnwidth]{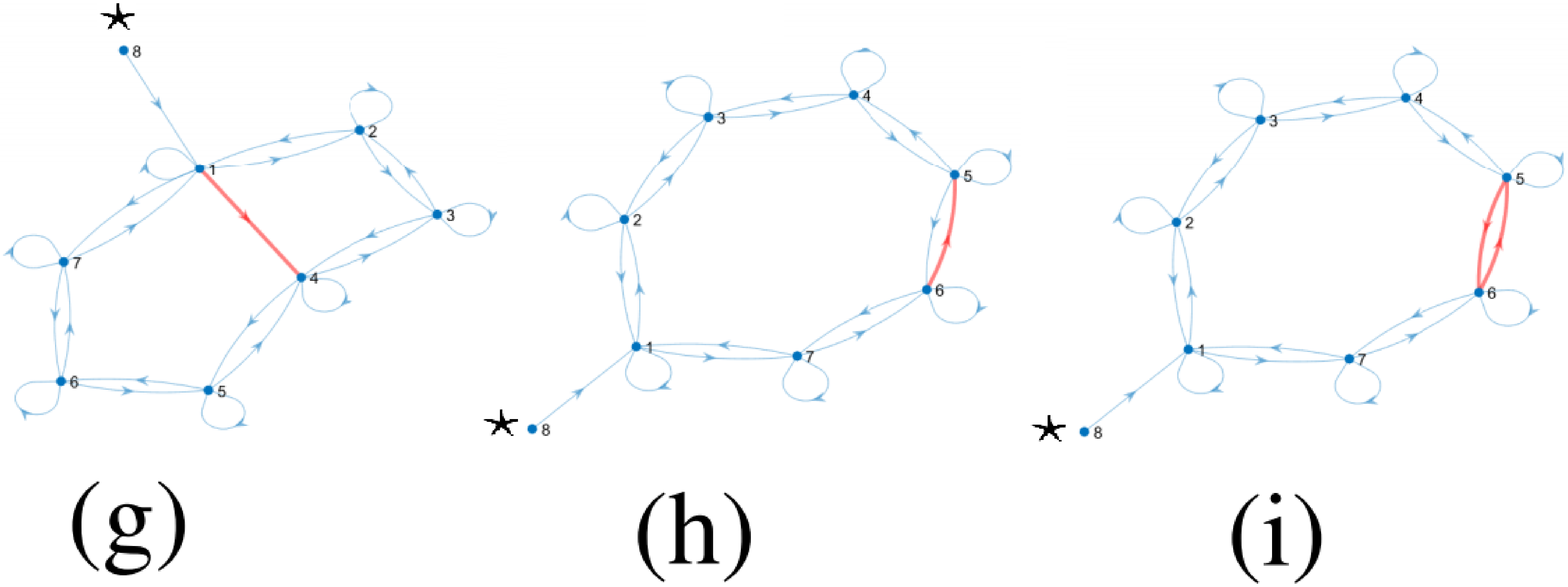}
	}

	\caption{The line, star, and circle networks in Example \ref{three-type}. Red edges are chosen to be perturbed.}  
	\label{example-three-typical}
\end{figure}

\end{example}

\begin{figure}
	\centering
	\includegraphics[width=3.5in]{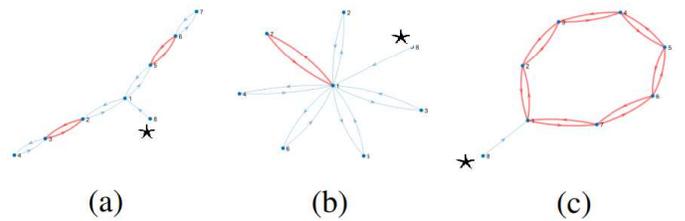}\\
	\caption{The networks in Example \ref{example-three-Laplacian}. Red edges are chosen to be perturbed.}\label{three-Laplacian}
\end{figure}

\begin{example} \label{example-three-Laplacian}
Finally, consider again the typical networks above with Laplacian dynamics; that is, the dynamics of the networks is $\dot x(t)= -L({\cal G}_{\star})x(t)+Bu(t)$, where $L({\cal G}_{\star})$ is the Laplacian matrix associated with ${\cal G}_{\star}$, $\star={\rm star, line, circ}$. Note for a network ${\cal G}_{\star}$ with $n$ nodes, $[L({\cal G}_\star)]_{ij}=-1$ if $(j,i)\in E({\cal G}_\star)$, $i\ne j$, and $[L({\cal G}_\star)]_{ii}=\sum \nolimits_{j=1,j\ne i}^n [L({\cal G}_\star)]_{ij}$. Due to the symmetry of ${\cal G}_{\star}$, it can be verified that all the $(-L(G_\star),B_\star)$ are uncontrollable for $\star={\rm star, line, circ}$. We use Algorithm \ref{alg2} to determine the set of edges whose weights need to be perturbed so that the resulting system becomes controllable, assuming symmetric edge weights (i.e., $[L(G)]_{ij}=[L(G)]_{ji}$ is always preserved). The parameters are the same as in Table \ref{parameters}. The evolution of $F(\theta^{(k)}, W^{(k)}, Z^{(k)})$ is given in Fig. \ref{convergence} (labelled `$\star$ Laplacian', $\star={\rm star, line, circ}$), and the obtained perturbed structures are presented in Fig. \ref{three-Laplacian}. It is worth mentioning that, although Algorithm \ref{alg2} is shown to converge by Fig. \ref{convergence}, the returned solutions in Fig. \ref{three-Laplacian} may not be the optimal ones.
\end{example}


\section{Conclusions}
This paper introduces and addresses the problem of determining the minimum number of entries in system matrices that need to be perturbed to make a given uncontrollable system controllable, i.e., the ZNCD.
It is shown computing the ZNCD is NP-hard, even when only the state matrices can be perturbed. Some lower and upper bounds of ZNCD are then given. Two heuristic algorithms are provided for computing ZNDC, the first of which is a greedy algorithm based on the structural controllability of a linearly parameterized plant, and the second one is built on the weighted $l_1$-norm relaxation and the sequential convex program. Those algorithms are valid when certain structural constraints are imposed on the corresponding perturbations. Finally, several numerical examples demonstrate the effectiveness of the proposed algorithms.

	%
	%
	
	%
	%
	%
	%
	%
	%
	
	
	{\footnotesize
		\bibliographystyle{elsarticle-num}
		\bibliography{yuanz3}
	}
\end{document}